%Paper: alg-geom/9404004
%From: ogrady@math.upenn.edu
%Date: Wed, 6 Apr 94 11:14:15 -0400

% 28 Marzo 1994

\def\a{\alpha}

\def\brel{\buildrel}
\def\bsk{\bigskip}

\def\CC{{\bf C}}

\def\cC{{\cal C}}

\def\cE{{\cal E}}
\def\cF{{\cal F}}
\def\cG{{\cal G}}

\def\cK{{\cal K}}
\def\cl{\colon}
\def\cL{{\cal L}}
\def\cM{{\cal M}}

\def\cO{{\cal O}}
\def\cod{{\rm cod}}
\def\cQ{{\cal Q}}
\def\cR{{\cal R}}
\def\cS{{\cal S}}

\def\d{\delta}
\def\D{\Delta}

\def\deg{{\rm deg}}

\def\det{{\rm det}}

\def\e{\epsilon}

\def\es{\emptyset}

\def\ker{{\rm ker}}

\def\G{\Gamma}

\def\hb{\hbox}
\def\hf{\hfill}

\def\hra{\hookrightarrow}

\def\i{\iota}
\def\id{{\rm id}}

\def\Im{{\rm Im}}

\def\ker{{\rm ker}}

\def\l{\lambda}

\def\L{\Lambda}
\def\lra{\longrightarrow}
\def\msk{\medskip}
\def\n{\noindent}

\def\o{\omega}
\def\ot{\otimes}
\def\op{\oplus}
\def\ov{\overline}
\def\O{\Omega}
\def\pf{\noindent{\bf Proof.}\hskip 2mm}
\def\PP{{\bf P}}

\def\qed{\hfill{\bf q.e.d.}}
\def\QQ{{\bf Q}}

\def\RR{{\bf R}}
\def\s{\sigma}
\def\S{\Sigma}

\def\ss{\subset}

\def\t{\theta}
\def\T{\Theta}

\def\tm{\times}

\def\ul{\underline}
\def\vf{\varphi}

\def\wh{\widehat}
\def\wt{\widetilde}

\def\ZZ{{\bf Z}}

%Fontenlargements.

\font\xxlbf=cmbx10 at 16pt

\font\xlbf=cmbx10 at 14pt

\font\lbf=cmbx10 at 12pt
\font\lrm=cmr10 at 12pt

\def\cMx{\cM_{\xi}}
\def\detx{\det_{\xi}}
\def\Dx{\D_{\xi}}
\def\cFx{\cF_{\xi}}
\def\cGx{{\cal G}_{\xi}}
\def\pal{\partial}
\def\Px{P_{\xi}}
\def\rx{r_{\xi}}
\def\Wx{W_{\xi}}

\magnification=1100

\centerline{\xxlbf Moduli of vector bundles on projective}
\msk
\centerline{\xxlbf surfaces: some basic results.}
\bsk
\bsk
\centerline{\lrm Kieran G. O'Grady}
\bsk
\centerline{\lrm University of Pennsylvania}
\msk
\centerline{\lrm Department of Mathematics}
\msk
\centerline{\lrm Philadelphia, PA 19104}
\msk
\centerline{\lrm e-mail: ogrady@math.upenn.edu}
\bsk
\centerline{\lrm March 28 1994}
\bsk
\bsk
Vector bundles on algebraic surfaces have been studied since
the 1960's.  Moduli spaces for stable bundles were
constructed in the 70's. Subsequently, Gieseker
and Maruyama constructed moduli spaces for
semistable torsion-free sheaves: these provide natural
compactifications of the moduli spaces of vector bundles. Many detailed
and interesting results have been proved regarding these moduli spaces
when the surface belongs to some particular class, for example if it is
$\PP^2$ or a $K3$, but  our
knowledge decreases as the Kodaira dimension of the surface increases,
and in particular very little is known  if the surface is of general type. In
this paper we address two basic questions:
\msk
\item{1.} When is the moduli space reduced and of the expected dimension?
\item{2.} When is the moduli space irreducible?
\msk
\n
In order to present our results we need to introduce some notation.
Let $S$ be  a smooth irreducible projective surface over $\CC$, and  let
$H$ be  an ample divisor on $S$.
To define a moduli space of sheaves
on the polarized surface $(S,H)$ we need a {\it set of sheaf data} $\xi$,
i.e.~a triple
$$\,\,\xi=(\rx,\detx,c_2(\xi))\,,$$
where $\rx$ is a positive integer,  $\detx$ is a line bundle on $S$, and
$c_2(\xi)\in
H^4(S;\ZZ)\cong\ZZ$. We let $\cMx$ be the moduli space of   semistable (with
respect to $H$) torsion-free sheaves, $F$, on $S$ with
$$\,\,r_F=r(F)=\rx\,,\quad\det F\cong\detx\,,\quad\
c_2(F)=c_2(\xi)\,.\eqno(0.1)$$
A fundamental theorem of Gieseker and Maruyama~[G1,Ma] asserts that
$\cMx$ is  projective.  If $F$ is a semistable sheaf
satisfying~(0.1), we let $[F]$ be the point in $\cMx$ corresponding to the
equivalence class of $F$.  We  recall some known facts concerning
the local structure of $\cMx$. First let's define the {\it discriminant} of  a
torsion-free sheaf $F$ on $S$ as
$$\D_F:=c_2(F)-{r(F)-1\over 2r(F)}c_1(F)^2$$
(warning: our normalization  differs from that  of~[DL]).
If $\xi$ is a set of sheaf data, the discriminant $\Dx$ is defined in the
obvious way: if $[F]\in\cMx$, then $\Dx=\D_F$.  The {\it expected dimension}
of $\cMx$ is given by
$$\,\,\hb{expdim}\left(\cMx\right):=2\rx\Dx-(\rx^2-1)\chi(\cO_S)\,.$$
Now assume that $[F]\in\cMx$ and that $F$ is  stable. Then
deformation theory~[F] gives
$$\eqalignno{\dim_{[F]}\cMx\ge &
2\rx\Dx-(\rx^2-1)\chi(\cO_S)\,,&(0.2)\cr
\dim T_{[F]}\left(\cMx\right) = &
 2\rx\Dx-(\rx^2-1)\chi(\cO_S)+h^0(F,F\ot K)^0\,,&(0.3)}$$
where, for a line bundle $L$ on $S$, we set
$$h^0(F,F\ot L)^0:=\dim\{\vf\in Hom(F,F\ot L)\ |\ {\rm tr}\vf=0\}\,.$$
Let $[F]\in\cMx$. Following Friedman we say that $\cMx$ is {\it good} at
$[F]$ if $F$ is stable and  $h^0(F,F\ot K)^0$ vanishes, where $K$ is the
canonical line bundle. In this case~(0.2) is an equality and the moduli
space is smooth near $[F]$.   We say that $\cMx$ is  {\it good} if it is good
at the generic point of every one of its irreducible components: this means
that   $\cMx$ is reduced and its dimension equals the expected dimension.
 Now we can go back
to Questions~(1) and~(2). First of all notice that if $\rx=1$ then $\cMx$ is
good (at each point) for trivial reasons, and furthermore, as is well-known,
it is always irreducible. Thus we will only be concerned with the case
$\rx\ge 2$. Our main result is that if $\Dx\gg 0$ then $\cMx$ is good and
irreducible, i.e.~both questions have a positive answer. The significance of
the condition $\Dx\gg 0$ is the following: if $\Dx<0$ then $\cMx$ is empty
by Bogomolov's Inequality, and on the other hand $\cMx\not=\es$ if $\Dx\gg
0$ (see~[HL,LQ]). Actually we will prove more than the simple statement
that $\cMx$ is asymptotically good. To explain this, let $L$ be a line bundle
on $S$ and  set   %
$$\,\,\Wx^L=\left\{[F]\in\cMx\ |\ h^0(F,F\ot L)^0>0\right\}\,.$$
Thus if $F$ is stable then $\cMx$ is good at $[F]$ if and only if $[F]\notin
\Wx^K$. We prove that (if $L$ is fixed) the growth of $\dim \Wx^L$ (for
fixed rank and  increasing $\Dx$) is smaller than that of the expected
dimension of $\cMx$. The theorem about $\cMx$ being good for $\Dx\gg 0$
follows at once from this result (setting $L=K$), together with some
dimension counts to take care of properly semistable sheaves. To see
that it is interesting to bound the dimension of $\Wx^L$  for  arbitrary
$L$  consider the case $L=\cO_S(K+C)$, where $C$
is a smooth curve on $S$.  In this case, if $F$ is locally-free and stable, the
geometric significance of
 $[F]\notin \Wx^L$ is the following:   the natural morphism from a
neighborhood of $[F]$ in $\cMx$ to the deformation space of $F|_C$
surjects onto the subspace of deformations fixing the isomorphism class
of $\det(F|_C)$.

The precise statements of the results we have described
are given in Theorems~B, C, D, E and their corollaries. One feature of
these theorems is that they are for the most part effective. Thus we give
an explicit upper bound for $\dim \Wx^L$. From this one can
 compute  an explicit lower bound for $\Dx$ guaranteeing that $\cMx$ is
good. This lower bound for arbitrary rank might not be very practical,
however it appears to depend on the "correct" quantities. If the
rank is two  our methods are somewhat stronger: in this
case we have computed the lower bound explicitely, and  we will show
that it can not be too far off from the optimal one. Regarding irreducibility
our results are less explicit, but we do give a lower bound for $\Dx$
guaranteeing irreducibilty of the moduli space for rank-two bundles with
trivial determinant on a complete intersection.

All of the above results spring from Propositions~(1.1)-(1.2). To
explain the content of these propositions let the
{\it boundary} of $\cMx$,  denoted by $\partial\cMx$, be the subset of
$\cMx$  parametrizing  sheaves which are singular, i.e.~not locally-free.
Furthermore,  if $X\ss\cMx$ let the {\it boundary} of $X$ be the
intersection
$$\,\,\partial X:=X\cap\partial\cMx\,.$$
Propositions~(1.1)-(1.2) assert that if $X\ss\cMx$ is a
closed subvariety whose dimension satisfies certain conditions then
$\pal X\not=\es$. These propositions are proved by  further developing the
ideas in the proof of Theorem~(1.0.3) of~[O]. (This theorem states that, in
rank two, any irreducible component intersects the boundary, if $\Dx\gg 0$.)
{}From Propositions~(1.1)-(1.2) one obtains Theorem~A, which
bounds the maximum dimension of complete subvarieties of $\cMx$ not
intersecting the boundary. Given Theorem~A one can easily bound
$\dim\Wx^L$, this is the content of Theorem~B, and in particular prove that
$\cMx$ is asymptotically good. Theorem~C gives an explicit lower
bound for $\Dx$ guaranteeing that $\cMx$ is good, in rank two: this  is
obtained by  arguments similar to those that give Theorem~B. Asymptotic
irreducibility (Theorem~D) follows from Theorems~A and~B: the argument
is due to Gieseker and Li. We reproduce their proof because we will then
apply it to complete intersections in order to obtain an explicit result
(Theorem~E).  Our proofs depend on certain estimates: in particular we
need an upper bound for the dimension of  the loci in $\cMx$ parametrizing
sheaves with subsheaves with (relatively) large slope. This and other
estimates are proved in the last section of the paper.

The first to prove that $\cMx$ is asymptotically good, if the rank is two,
was Donaldson~[D] (see also~[F,Z]). He proved that $\dim\Wx^K$ is bounded,
up to lower order terms, by $3\Dx$ (the same proof works also  for
$\Wx^L$);  by~(0.2) this implies that $\cMx$ is good for $\Dx\gg 0$.
Donaldson's bound for $\dim\Wx^L$ is asymptotically better than that given
by Theorem~B, but since the lower order terms in his formula have eluded
computation, it does not give an effective result. Recently Gieseker and
Li~[GL2] have proved that $\cMx$ is asymptotically good. They proved
that the codimension of $\Wx^L$ in $\cMx$ goes to infinity for $\Dx\gg 0$,
but they did not make their result  effective. Finally Gieseker and Li~[GL1]
were the first to prove  asymptotic irreducibility in rank two.

\bsk
\n
{\bf Statement of results.}
\hskip 2mm
Throughout the paper surface means
a smooth irreducible projective surface: we will always denote it by $S$.
We let $K$ be its  canonical divisor, and   $H$ be  an ample divisor on
$S$. We will often make the following  assumption:
$$\,\,\hb{$|H|$ is base-point-free and $\dim|H|\ge 2$.}
\eqno(0.4)\,.$$
When considering a moduli space $\cMx$ we always tacitly assume that
$\rx\ge 2$. If $r>2$ is an integer, set
$$\eqalign{\rho(r) &:=8(16r^3-39r^2+36r-12)^{-1}\,,\cr
\D_0(r,S,H)&:=\rho^{-1}H^2\,,\cr
\l_2(r)&:= 2r-{r-1\over 2}\rho\,,\cr
\l_1(r,S,H)& := \sqrt{\rho}\left[{r^3\over 2}{|K\cdot H|\over\sqrt{H^2}}+
r^5\sqrt{H^2}\right]\,,\cr
\l_0(r,S,H)&:={r^7\over 2}H^2+{r^4\over 5}{(K\cdot H)^2\over H^2}
+{r^2\over 2}{(K\cdot H+(r^2+1)H^2+1)^2\over H^2}
+r^2|\chi(\cO_S)|+{r^3\over 8}|K^2|\,.}$$
When $r=2$, we set
$$\eqalign{\D_0(2,S,H):= & \cases{ 3H^2\,, & if $K\cdot H<0$, \cr
3H^2\left(1+{K\cdot H\over H^2}\right)^2\,, & if $K\cdot H\ge 0$,}\cr
\l_2(2):= & {23\over 6}\,,\cr
\l_1(2,S,H):= & {1\over 2\sqrt{3H^2}}(4H^2+3K\cdot H+4)\,,\cr
\l_0(2,S,H):=&\cases{{3(K\cdot H+H^2+1)^2\over 2H^2}+{(K\cdot H)^2\over
4H^2}-{K^2\over 4}+4-3\chi(\cO_S)\,, & if $K\cdot H<0$,\cr
{3(K\cdot H)^2\over H^2}+6K\cdot H+{3H^2\over 2}-{K^2\over4}
+8-3\chi(\cO_S)\,, & if $K\cdot H\ge 0$.}}$$

\proclaim Theorem A.
Let $(S,H)$ be a polarized  surface. Assume
that $H$ satisfies~(0.4).   Let $\xi$ be a set of sheaf data such that
$\D_{\xi}>\D_0(\rx,S,H)$.  If $X\ss\cM_{\xi}$ is a closed
subvariety such that
$$\,\,\dim X>\l_2(\rx)\D_{\xi}+\l_1(\rx,S,H)\sqrt{\D_{\xi}}
+\l_0(\rx,S,H)\,,\eqno(0.5)$$
then the boundary of $X$ is non-empty.

\n
Notice that the above result is meaningful because $\l_2(r)<2r$
(see~(0.2)).
Let
$$\eqalign{\l_0'(r,S,H):= & \max\left\{\l_0(r,S,H),\e(r,S,H)+r\right\}\,,\cr
\D_1(r,S,H):=&\max\left\{\left(\l_2-(2r-1)\right)^{-1}\cdot
\left(\D_0+e_K-\l_0'-(r^2-1)\chi(\cO_S)\right),\D_0\right\}\cr
& \hb{if $(r,S,H)\not=(2,\PP^2,\cO_{\PP^2}(1))$,}\cr
\D_1(2,\PP^2,\cO_{\PP^2}(1)):=& 3\,.}$$
Here $\e$ and $e_K$ are given by~(5.39) and~(5.6) respectively.

\proclaim Theorem B.
Let $(S,H)$ be a polarized surface, with $H$ satisfying~(0.4).   Let   $L$ be
a line bundle on $S$. Let $\xi$ be a set of sheaf data such that
$\Dx>\D_1(\rx,S,H)$. Then
$$\,\,\dim\Wx^L\le
\l_2(\rx)\Dx+\l_1(\rx,S,H)\sqrt{\Dx}+\l_0'(\rx,S,H)+e_L(\rx,S,H)\,,
\eqno(0.6)$$
where  $e_L$ is given by~(5.6).

\n
Applying the above result with $L=K$, one gets the following

\proclaim Corollary ${\rm B}'$.
There exists a function $\D_1'(r,S,H)$ (depending only on $r$,
$K^2$, $K\cdot H$, $H^2$,  and $\chi(\cO_S)$) such that the following holds.
Let $(S,H)$ be as above, and $\xi$ be a set of sheaf data such that
$\Dx>\D_1'(\rx,S,H)$. Then $\cMx$ is good. Furthermore $\cMx$ is the closure
of the open subset parametrizing  $\mu$-stable vector bundles.

\n
In the rank-two case we will carry out the computations
necessary to determine explicitely the lower bound of this
corollary. The result is the following

\proclaim Theorem C.
Let $(S,H)$ be a polarized  surface. Assume that $H$ is
effective and that there exists a smooth curve in $|n_0H|$, where  $n_0$
 is given by~(2.13). Let $\D_2(S,H)$ be the
function defined by~(2.12). Let $\xi$ be a set of
sheaf data with $\rx=2$. If
$$\,\,\Dx>\D_2(S,H)+2h^0(2K)\,,$$
then $\cMx$ is good. Furthermore $\cMx$ is the closure of the open subset
parametrizing  $\mu$-stable vector bundles.

\n
To exemplify the possible uses of this theorem we give
the following  two results.

\proclaim Corollary~${\rm C}'$.
Let $S$ be a surface with ample canonical divisor  $K$. Assume that
$p_g(S)>0$, and that $K^2\gg 0$ ($K^2>100$ will do). Let $\xi$ be a set of
sheaf data with $\rx=2$, and let $\cMx$ be the corresponding moduli space
for semistable sheaves on the polarized surface $(S,K)$. If
$$\,\,\Dx\ge 42K^2+15\chi(\cO_S)\,,$$
then $\cMx$ is good, and furthermore  the subset parametrizing
$\mu$-stable vector bundles is dense in $\cMx$.

\proclaim Corollary~${\rm C}''$.
Let $S$ be a surface with ample canonical divisor  $K$. Assume that $H$ is
a divisor satisfying~(0.4), and that $c_1(K)=kc_1(H)$ for $k\gg 0$
(say $k>100$). Let $\xi$ be a set of sheaf data with $\rx=2$, and let
$\cMx$ be the corresponding moduli space for semistable sheaves on
$(S,K)$. If
$$\,\,\Dx\ge17K^2+10\chi(\cO_S)\,,$$
then the same conclusions as in the previous corollary hold.

\n
These two corollaries only apply to minimal surfaces of general type with
no $(-2)$-curves. Of course Theorem~C applies to any surface, in particular
we can apply it to general-type surfaces with $(-2)$-curves, but the
lower bound one gets is not  as nice as that of Corollaries~${\rm C}'$-${\rm
C}''$. However we believe that bounds similar to those of
these corollaries  hold for any surface of general type,
if $K^2$ is replaced by $\o^2_{S_{can}}$, where $S_{can}$ is the canonical
model of $S$, and the polarization is close enough to $\o_{S_{can}}$ (how
"close "will depend on $\Dx$). When $S_{can}$ is smooth this should
follow from the results in Section~2 of~[MO]; the case when $S$ contains
$(-2)$-curves should be analyzable by similar methods.

Finally we come to irreducibility.

\proclaim Theorem~D.
There exists a function $\D_3(r,S,H)$ such that the following
holds. If $(S,H)$ is a polarized surface, and if $\xi$ is set of
sheaf data such that $\Dx>\D_3(\rx,S,H)$, then $\cMx$ is
irreducible (and is the closure of the locus parametrizing
$\mu$-stable vector bundles).

\n
We give an explicit value for $\D_3$ valid for complete
intersections of large degree, if the rank is two and the determinant is
trivial.

\proclaim Theorem~E.
Let $S$ be a complete intersection in a projective space, and let
$H=\cO_S(1)$. Suppose also that the integer $k$ such that
$K\sim kH$ is very large ($k>100$ suffices). Let
$$\,\,\xi=\left(2,\cO_S,c_2(\xi)\right)\,.$$
If
$$\,\,\Dx>95K^2+11\chi(\cO_S)+1\,,$$
then $\cMx$ is irreducible.

\bsk

\n
{\bf Notation and conventions.}
\msk
\n
We work throughout over the complex numbers. Sheaves are always
coherent.  We let $\pi_X\cl X\tm Y\to X$ be  the projection. A {\it family
of sheaves on $X$ parametrized by $B$} consists of a sheaf $\cF$  on
$X\tm B$, flat over $B$.   If $b\in B$ we set $\cF_b:=\cF|_{X\tm \{b\}}$.
Let $\xi$ be a set of sheaf data for $S$, and let $U\ss\cMx$ be an
algebraic  subset   all of whose points parametrize stable sheaves. A {\it
tautological family}  parametrized by $U$ consists of a family of stable
sheaves $\cF$ on $S$, parametrized by $U$, such that for all $b\in U$ the
isomorphism class of $\cF_b$ is represented by $b$. If $F$ is a sheaf on a
smooth curve or surface, we let $Def(F)$ be the versal deformation space
of $F$, and $Def^0(F)$ be the subscheme parametrizing deformations which
fix the isomorphism class of $\det F$. (See~[F].)

If $F$ is a semistable torsion-free sheaf on $S$, we denote by $Gr(F)$  the
direct sum of the successive quotients of any Jordan-H\"older filtration
of $F$. Recall that the  (closed) points of $\cMx$ are in one-to-one
correspondence with equivalence classes of torsion-free semistable
sheaves, $F$, satisfying~(0.1): two semistable sheaves $F_1$, $F_2$
are  equivalent if  $Gr(F_1)\cong Gr(F_2)$.

Let $X$ be a projective irreducible variety,
and $D$ be an ample divisor   on $X$. The {\it slope} of a torsion-free sheaf
$F$ on $X$ with respect to $D$ is given by
$$\,\,\mu_F=\mu(F):={1\over r_F}c_1(F)\cdot D^{n-1}\,,$$
where $n:=\dim X$. We recall that $F$ is $\mu$-{\it semistable}
(equivalently {\it slope-semistable}) if, for all subsheves $E\ss F$ we have
$$\,\,\mu_E\le\mu_F\,.$$
If the above inequality is strict whenever $r(E)<r(F)$, then $F$ is
$\mu$-{\it stable} ({\it slope-stable}). If $F$ is semistable
then it is $\mu$-semistable, and  $\mu$-stability implies stability.
Let $\a\in\RR$. A torsion-free sheaf $F$ on $S$ is {\it $\a$-stable} if, for
every subsheaf $E\ss F$ with $0<r(E)<r(F)$, one has
$$\,\,\mu_E<\mu_F-{\a\over r_{\cE}}\sqrt{H^2}\,.$$
Thus  $\mu$-stability is equivalent to
$0$-stability. As is immediately verified $F$ is $\a$-stable if and only if,
for every non-trivial torsion-free quotient $F\to Q$, one has
$$\,\,\mu_Q>\mu_F+{\a\over r(Q)}\sqrt{H^2}\,.$$
Furthermore the notion of $\a$-stability only depends on the ray spanned by
$c_1(H)$, and  $F$ is $\a$-stable if and only if so is $F^{*}$. We let
$$\,\,\cMx(\a):=\{[F]\in\cMx\ |\ \hb{$F$ is $\ul{\rm not}$
$\a$-stable}\}\,.$$
We will show (Proposition~(5.9)) that $\cMx(\a)$ is a constructible
subset of $\cMx$. In case $\rx=2$ we will also consider a
(constructible) subset of $\cMx(\a)$, defined as follows.

\proclaim (0.7) Definition.
Let $\xi$ be a set of sheaf data with $\rx=2$.  Let  $C\ss
S$ be a smooth irreducible curve. For $\a\in\RR$ we  let
$\cMx^C(\a)\ss\cMx$ be the subset of points $[F]$ such that $F|_C$ is
locally-free, and such that there exists a rank-one subsheaf  $A\ss F$ with:
\msk
\item{1.} $\mu_A\ge \mu_F-\a\sqrt{H^2}$ (so $F$ is not $\a$-stable).
\item{2.} The restriction $A|_C$ spans a  destabilizing subline bundle of
$F|_C$.

\bsk
\bsk

\n
{\xlbf 1. A criterion for non-emptiness of the boundary.}
\msk
\n
In this section we will  prove  the two following propositions.

\proclaim (1.1) Proposition.
Let $(S,H)$ be a polarized surface, with $H$ satisfiying~(0.4).  Let $\xi$ be
a set of sheaf data, and let $X\ss\cMx$ be a closed  subvariety. Assume that
there exists a positive integer $n$ such that:
\msk
\item{1.} $\dim X>{1\over 2}(\rx^2-1)(H^2n^2+K\cdot Hn)$,
\msk
\item{2.} $\dim X>{1\over 8}\rx^2(H^2n^2+K\cdot Hn)+
{1\over 4}\rx^2+{1\over 2}+(2\rx-1)\Dx+\e(\rx,S,H)$,
\msk
\item{3.} $\dim X>{1\over 8}\rx^2(H^2n^2+K\cdot Hn)+
{1\over 4}\rx^2+{1\over2}+
\dim\cMx\left((\rx-1)\sqrt{H^2}n\right)$,
\msk
\item{4.} $\dim
X>2\rx\Dx-(r_{\xi}^2-1)\chi(\cO_S) +e_K(\rx,S,H) +{1\over
4}r_{\xi}^2-{1\over 2}(r_{\xi}-1)\left(H^2n^2-K\cdot Hn\right)$,
\msk
\n
where $[*]$ denotes integer part, and $e_K(\rx,S,H)$, $\e(\rx,S,H)$  are as
in~(5.6) and~(5.39) respectively. Then  $\partial X\not=\es$.

\proclaim (1.2) Proposition.
Let $(S,H)$ be a polarized surface. Assume that $H$ is effective.
Let $\xi$ be a set of sheaf data with $\rx=2$.   Let
$X\ss\cMx$ be a closed  subvariety. Suppose that there exist a positive
integer $n$ and a smooth curve $C\in |nH|$ such that Items~(1) and~(2)  of
Proposition~(1.1) are satisfied, and  furthermore
$$\eqalignno{\dim X> & {1\over 2}(H^2n^2+K\cdot Hn)+1+
\dim\cMx^C\left(\sqrt{H^2}n\right)\,,&(1.3)\cr
\dim X> & {1\over 2}(H^2n^2+K\cdot Hn)+1+
\dim\cMx\left({K\cdot H\over 2\sqrt{H^2}}\right)\,,&(1.4)\cr
\dim X> & 4\Dx-3\chi(\cO_S) +h^0(2K)+1
-{1\over 2}\left(H^2n^2-K\cdot Hn\right)\,.&(1.5)}$$
Then  $\partial X\not=\es$.

\n
The above propositions will be proved at the end of the section.
\bsk

\n
{\bf Certain families of elementary modifications.}
\msk
\n
In this subsection we will prove:

\proclaim (1.6) Proposition.
Let $(S,H)$ be a polarized surface. Let $\xi$ be a set of sheaf data.
Let $[F]\in\cM_{\xi}$, and assume that $F$ is locally-free and $\mu$-stable.
Let $C\ss S$ be a smooth irreducible curve. Set $\a_C:=(\rx-1)(C\cdot
H)/\sqrt{H^2}$. Assume that:
\msk
\item{1.}  $F|_C$ is  not stable,
\item{2.} $[F]\notin\cMx\left(\a_C\right)$.
\msk
\n
Let $X\ss\cM_{\xi}$ be a closed  subvariety containing $[F]$, and
such that
$$\,\,\dim X>2r_{\xi}\D_{\xi}-(r_{\xi}^2-1)\chi(\cO_S)+h^0(F,F\ot K_S)^0
+{r_{\xi}^2\over 4}-{1\over 2}(\rx-1)C^2
+{1\over 2}(\rx-1)C\cdot K\,.\eqno(1.7)$$
Then the boundary of $X$ is non-empty. If $\rx=2$, the same conclusion holds
if Item(2) is replaced by the condition
$$\,\,[F]\notin\cMx^C\left(\a_C\right)\,.$$

We will prove this proposition   at the end of the subsection. The key
ingredient in its proof is  provided by a certain family of elementary
modifications. We now proceed to introduce this family.

Let   $C\ss S$ be as above. Let $[F]\in\cMx$, and
assume that $F|_C$ satisfies Item~(1) of Proposition~(1.6).
Choose a destabilizing sequence
$$\,\,0\to \cL_0\to F|_C\to \cQ_0\to 0\,.\eqno(1.8)$$
By definition we have:
$$\displaylines{\hf\hb{$\cL_0$ and $\cQ_0$ are locally-free,}\hf\llap{(1.9)}\cr
\hf \,\,\mu(\cL_0)-\mu(\cQ_0)\ge 0\,.\hf\llap{(1.10)}}$$
Let $E$ be the elementary modification of $F$ associated to the
destabilizing quotient of~(1.8), i.e.~the sheaf on $S$ fitting into  the exact
sequence
$$\,\,0\to E\to F\brel g\over \lra\i_*\cQ_0\to 0\,,\eqno(1.11)$$
where $\i\cl C\hra S$ is the inclusion.  Restricting the above sequence to $C$,
one gets
an exact sequence
$$\,\,0\to\cQ_0\ot\cO_C(-C)\to
 E|_C\brel f_0\over\lra\cL_0\to 0\,.\eqno(1.12)$$
Hence by~(1.9) $E|_C$ is locally-free. Since $E$ and $F$ are isomorphic
outside of $C$ we conclude that $E$ is locally-free.  Let
$$Y_F:=Quot(E|_C;\cL_0)$$
 be the Grothendieck $Quot$-scheme parametrizing quotients of $E|_C$ which
have the same Hilbert polynomial as $\cL_0$.  Notice that the notation is
slightly imprecise, since a destabilizing sequence for $F|_C$ is not
necessarily unique. However this   will not create confusion because we will
always fix a sequence~(1.8) once and for all. We  denote  by $0$ the point
of $Y_F$ corresponding to  $f_0$ (see~(1.12)). We are now ready
to define the promised family of elementary modifications. The parameter
space will be $Y_F$. Let
$$\pi_C^*(E|_C)\brel f\over\to\cL$$
 be the tautological quotient sheaf on $C\tm Y_F$, and   let $\cG$ be the sheaf
on $S\tm
Y_F$ fitting into the exact sequence
$$\,\,0\to \cG\to\pi_S^*E\brel \phi\over\to (\i\tm \id_{Y_F})_*\cL\to 0\,,$$
where $\phi$ is the composition of restriction to $C\tm Y_F$ and $f$.

\proclaim (1.13) Lemma.
The sheaf $\cG$ is flat over $Y_F$, and thus we can regard it as a family of
sheaves on $S$ parametrized by $Y_F$. Let $y\in Y_F$. Then $\cG_y$ fits
into the exact sequence  %
$$\,\,0\to \cG_y\to E\brel \phi_y\over\lra \i_*\cL_y\to 0\,,\eqno(1.14)$$
where $\phi_y:=\phi|_{S\tm\{y\}}$. In particular $\cG_y$ is torsion-free.
Finally, $\cG_y$ is
singular if and only if so is $\cL_y$.

\pf
The sheaf $\cG$ is $Y_F$-flat because so are $\pi_S^*E$ (obvious) and
$(\i\tm\id_{Y_F})_{*}\cL$ (by definition of the $Quot$-scheme). Flatness of the
latter
implies that~(1.14) is exact. For $y\in Y_F$ let $\cR_y:=\ker f_y$; since
$E|_C$ is locally-free, so is $\cR_y$.  The exact sequence
$$0\to \cL_y(-C)\to \cG_y|_C\to \cR_y\to 0$$
shows that $\cG_y$ is singular at a point  $P\in C$ if and only if $\cL_y$ is
singular at $P$. Since
$$\,\,\cG_y|_{(S-C)}\cong E|_{(S-C)}\,,$$
 and $E$ is locally-free, we conclude that $\cG_y$ is singular if and only if
so is $\cL_y$.
\qed

\msk
Let $\cF:=\cG\ot\pi_S^*\cO_S(C)$. By the above lemma,   we can regard $\cF$ as
a
family of torsion-free sheaves on $S$ parametrized by $Y_F$. Let $\partial
Y_F\ss Y_F$ be the subset  parametrizing singular sheaves.

\proclaim (1.15) Lemma.
Let notation be as above. Then:
\item{1.} $\cF_0\cong F$.
\item{2.} For $y\in Y_F$ we have $\det\cF_y\cong\detx$ and
$c_2(\cF_y)=c_2(\xi)$.
\item{3.} Let $\S\ss Y_F$ be a closed  subvariety. If
$\dim \S>\rx^2/4$, then $\S\cap\partial Y_F\not=\es$.

\pf
 Clearly the subsheaf $F(-C)\hra F$ is in the kernel of the map $g$ of~(1.11).
Hence $F(-C)$
is actually a subsheaf of $E$; let $\l\cl F(-C)\to E$ be the inclusion map. As
is easily
checked
$$\,\,\Im\left(\l|_C\right)=\ker f_0\,.$$
Since $\l$ is an isomorphism outside of $C$, we conclude that $F(-C)$ fits into
the
exact sequence
$$\,\,0\to F(-C)\to E \brel \phi_0\over\lra \i_*\cL_0\to 0\,.$$
By Lemma~(1.13) the sheaf $\cG_0$ fits into the same  exact sequence, and thus
$\cG_0\cong F(-C)$. This proves Item~(1). Let's consider the second Item. It
follows from
Exact sequences~(1.11) and~(1.14) that $\det E\cong\det
F(-r_{\cQ_0}C)$ and that $\det \cG_y\cong \det E(-r_{\cL_y}C)$. This
gives $\det\cF_y\cong\det F\cong\detx$. Since $c(\i_*\cL_y)$ is
independent of $y\in Y_F$, so is $c(\cG_y)$, and hence also $c(\cF_y)$.
Since  $c_2(\cF_0)=c_2(F)=c_2(\xi)$, we conclude that
$c_2(\cF_y)=c_2(\xi)$ for all $y\in Y_F$. Now let's prove Item~(3). Assume
that $\S\cap\partial Y_F=\es$; we will arrive at a contradiction. Fix $P\in
C$. Since we are assuming that $\cF_y$ is locally-free for all $y\in\S$, it
follows from the last sentence of  Lemma~(1.13) that
$$\cL|_{P\tm\S}$$
is a vector bundle. Thus the kernel of the map
$$\,\,f|_{P\tm\S}\cl E_P\ot\cO_{\S}\to\cL|_{P\tm\S}\,,$$
is a rank-$r(\cQ_0)$ $\ul{\rm subbundle}$ of the trivial bundle on $\S$
with fiber $E_P$ (here $E_P$ is the  fiber of $E$
at $P$).   Let
$$\rho\cl\S\to {\bf Gr}:={\bf Gr}(r(\cQ_0),E_P)$$
be the morphism defined by setting $\rho(y):=\ker
\left(f_{(P,y)}\right)$. Now fix $[V]\in {\bf Gr}$.  By  hypothesis
$\dim\S>\dim {\bf Gr}$, and hence
$$\,\,\dim\rho^{-1}\left([V]\right)\ge 1\,.\eqno(*)$$
Set
$$\,\,\O:=
\bigcup_{y\in\rho^{-1}\left([V]\right)}\PP\left(\ker f_y\right)\,.$$
Since the map $\S\to\PP(E)$ defined by
$y\mapsto\PP\left(\ker f_y\right)$ is injective, we conclude by~($*$) that
$$\,\,\dim\O\ge (1+r(\cQ_0))\,.\eqno(\dag)$$
 Since $Y_F$ is complete so is $\S$, and hence   $\O$
 is a closed subvariety of $\PP(E)$.   By~($\dag$) we conclude that
$$\,\,\dim\O\cap\PP(E_P)\ge r(\cQ_0)\,.$$
This is absurd because the above intersection is $\PP(V)$, and
$\dim\PP(V)=(r_{\cQ_0}-1)$.
\qed

\msk
In order to use the above lemma, we need to ensure that the dimension
of $Y_F$ is large, and that $\cF$ is a family of  semistable sheaves. (Notice
that by Item~(1) the sheaves $\cF_y$ are stable for $y$ varying in an open
non-empty subset of $Y_F$; however this will not be good enough.)

\proclaim (1.16) Lemma.
Keep notation as above. Then
$$\,\,\dim Y_F\ge{1\over 2}(\rx-1)C^2-{1\over 2}(\rx-1)C\cdot K\,.$$

\pf
By~(1.12) we have the  lower bound
$$\,\,\dim Y_F\ge \chi\left(\cQ_0^{*}\ot\cO_C(C)\ot \cL_0\right)\,.$$
Let  $g$ be the genus of $C$. Riemann-Roch gives
$$\,\,\chi\left(\cQ_0^{*}\ot\cO_C(C)\ot \cL_0\right)
=r_{\cL_0}r_{\cQ_0}\left[\mu(\cO_C(C))+\mu(\cL_0)-\mu(\cQ_0)+1-g\right]\,,$$
where the slopes are as $\ul{\hb{bundles on $C$}}$. By
Inequality~(1.10) we conclude that
$$\,\,\dim Y_F\ge r_{\cL_0}r_{\cQ_0}\left(C^2+1-g\right)\,.$$
Using adjunction one gets the lemma. (Notice that if $(C^2-C\cdot K)<0$ then
the
lemma is trivially verified.)
\qed

\msk
\n
Regarding stability we have the following

\proclaim (1.17) Lemma.
Keep notation as above, and let $\a_C$ be as in the statement of
Proposition~(1.6). If $[F]\notin \cMx\left(\a_C\right)$, then $\cF$ is a
family of stable sheaves. In the case $\rx=2$ the same conclusion holds if:
\msk
\item{1.} $F$ is $\mu$-stable, and
\item{2.} $[F]\notin\cMx^C(\a_C)$.

\pf
Let $y\in Y_F$. We will show that $\cG_y$ is $\mu$-stable; this will prove the
lemma. First notice that, by Item~(2) of Lemma~(1.15), we have
$$\,\,\mu(\cG_y)=\mu_F-C\cdot H\,.\eqno(*)$$
Now let $A\hra\cG_y$ be a subsheaf with $0<r(A)<r(\cG_y)$. Let $\l\cl A\to
F$ be the composition (see~(1.14) and~(1.11))
$$\,\,A\hra\cG_y\to E\to F\,.$$
Since $\l$ is injective outside of $C$, and since $A$ is torsion-free, we
conclude that
$\l$ is an injection. If $[F]\notin\cMx\left(\a_C\right)$ then
$$\,\,\mu_A<\mu_F-{\a_C\over r_A}\sqrt{H^2}
\le \mu_F-C\cdot H=\mu(\cG_y)\,,\eqno(\dag)$$
and hence $\cG_y$ is $\mu$-stable. Now assume $\rx=2$ and
$[F]\notin\cMx^C(\a_C)$. If $\l$ is zero at the generic point of $C$, then we
get  an injection $A(C)\hra F$. By hypothesis $F$ is $\mu$-stable, and hence
$$\,\,\mu_A+C\cdot H<\mu_F\,.$$
By~($*$) we get $\mu_A<\mu(\cG_y)$. Now assume $\l$ is not zero at the
generic point of $C$. Then
$$\,\,\Im\left(\l|_C\right)=\Im\left(E|_C\to F|_C\right)
\quad\hb{at the generic point of $C$.}$$
Since the  right-hand side  is a destabilizing subline bundle of $F|_C$,
and since $[F]\notin\cMx^C(\a_C)$ we conclude that~($\dag$) holds. Thus
$\cG_y$ is $\mu$-stable.
\qed

\bsk

\n
{\it Proof of Proposition~(1.6).}
\hskip 2mm
Since $F$ satisfies Item~(1), we can construct $Y_F$ and $\cF$ as above.
By  Lemma~(1.17), $\cF$ is a family of stable
sheaves on $S$  parametrized by $Y_F$. Hence, by Item~(2) of Lemma~(1.15),
$\cF$ induces a classifying morphism
$$\,\,\vf\cl Y_F\to \cMx\,.$$
By Item~(1) of~(1.15), we have $\cF_0\cong F$. Hence $[F]\in X$ and   the
inverse image $\vf^{-1}X$ is a closed subvariety of $Y_F$ containing the
point $0$. We have
$$\,\,\dim\vf^{-1}X\ge \dim Y_F-\left(\dim T_{[F]}\cMx-\dim X\right)\,.$$
To be precise: the right-hand side is a lower bound for the dimensions of all
irreducible components of  $\vf^{-1}X$ containing $0$. By~(0.3), (1.7) and
Lemma~(1.16)  we conclude that $\dim\vf^{-1}X>\rx^2/4$. By Item~(3) of
Lemma~(1.15) there exists $y\in\vf^{-1}X$ such that $\cF_y$ is singular.
Then $\vf(y)\in\partial X$, and hence $\partial X\not=0$.
\qed

\bsk

\n
{\bf Determinant bundles.}
\msk
\n
In this subsection $C$ will be a smooth irreducible curve in the linear system
$|nH|$.
Let  $\cM(C;\xi)$ be the moduli space of
rank-$\rx$ semistable bundles on $C$ with determinant isomorphic to $\detx|_C$.
Let $X\ss \cMx$ be a subvariety such that, for all $[F]\in X$, the restriction
$F|_C$ is a
stable locally-free bundle. Since $C\in |nH|$ this implies that $F$ is
$\mu$-stable for all $[F]\in X$. Then, as is easily verified there exists a
morphism
$$\rho\cl X\to \cM(C;\xi)$$
given by restriction, i.e.~$\rho([F])=[F|_C]$. Our goal will be  to
prove the following

\proclaim (1.18) Proposition.
Let $X,C,\rho$ be as above. Assume also  that $X$ is closed and irreducible,
and that all sheaves parametrized by points of $X$ are locally-free.  Let
$\T$ be the theta-divisor on $\cM(C;\xi)$ (see~DN]). Then
$$\,\,\left(\rho^*\T\right)^{\dim X}>0\,.$$

We will first prove a series of preliminary results.  We start
with  a weak version of closedness of
non-stability for a family of vector bundles on a variable degenerating curve.
More precisely:
Let $B$ be a smooth curve, $0\in B$ be a base point, and $B^0:=(B-0)$. Let
$\pi\cl\cC\to B$ be
a family of  curves; for $b\in B$ set $C_b:=\pi^{-1}(b)$. We assume that:
\msk
\item{1.} $\cC$ is smooth outside a finite set of points in $\cC_0$,
\item{2.}  all fibers $C_b$ are reduced and connected, and  for $b\not=0$
they are smooth.
\msk
\n
Let $D_1,\ldots,D_s$ be the irreducible components of the
central fiber $C_0$.

\proclaim (1.19) Lemma.
Keep notation as above. Let $\cF$ be a vector bundle on $\cC$, and set
$\cF_b:=\cF|_{C_b}$. Assume that for all $b\not=0$ the bundle $\cF_b$ is
non-stable . There exists $1\le i\le s$ such that, letting $\l_i\cl
\wt{D}_i\to D_i$ be the normalization, the bundle
$\l_i^*\left(\cF|_{D_i}\right)$ is not stable.

\pf
For generic $b\in B^0$ the Harder-Narasimhan filtration of $\cF_b$ has constant
type
(i.e.~length, and rank and slope of the successive quotients). Thus, shrinking
$B^0$ if
necessary, we can assume that there exists a vector bundle $\cQ^0$ on
$\cC^0:=\pi^{-1}(B^0)$, and  an exact sequence
$$\,\,\cF|_{\cC^0}\brel\a\over\lra \cQ^0\to 0\,,\eqno(*)$$
whose restriction to $C_b$ is a destabilizing sequence, for all $b\in B^0$.
By properness of the relative $Quot$-scheme parametrizing quotients of $\cF_b$,
there is a
$B$-flat sheaf $\cQ$ on $\cC$  extending  $\cQ^0$, and an exact sequence
$$\,\,\cF\to\cQ\to 0\,.\eqno(\star)$$
By flatness $\cQ$ is torsion-free. In particular it is locally-free
outside a finite set of points in $C_0$.
Let $f\cl\wt{\cC}\to\cC$ be a desingularization such that, for all $1\le i\le
s$, the
proper tranform of $D_i$ is smooth.  We denote this proper transform by
$\wt{D}_i$.
Let $Tor(f^*\cQ)$ be the torsion subsehaf of $f^*\cQ$: since $\cQ$ is
torsion-free  $Tor(f^*\cQ)$  is supported on the exceptional divisors of $f$.
Pulling
back~($\star$) we get an exact sequence
$$\,\,0\to \cK\to f^*\cF\to f^*\cQ/Tor(f^*\cQ)\to 0\,.\eqno(\dag)$$
Since $\wt{\cC}$ is a smooth surface, and since $f^*\cF$ is locally-free,
$\cK$ is also
locally-free. Let $\vf:=\pi\circ f$, and set $\wt{\cC}_b:=\vf^{-1}(b)$. If
$b\not=0$, then
the restriction of~($\dag$) to $\wt{\cC}_b=\cC_b$ is the destabilizing
sequence associated to $\a_b$ (see~($*$)), and hence
$$\,\,{1\over r(\cK)}c_1(\cK)\cdot \wt{\cC}_b
\ge {1\over r(\cF)}c_1(f^*\cF)\cdot \wt{\cC}_b\,.\eqno(\sharp)$$
The same inequality holds also when $b=0$. We have
$$\,\,\wt{\cC}_0=\sum_{i=1}^s\wt{D}_i+\sum_{j=1}^pn_jE_j\,,\eqno(\flat)$$
where $E_1,\ldots,E_p$ are the exceptional divisors of $f$,  and $n_j>0$ for
all $j$. The map $\cK\to f^*\cF$ has isolated zeroes, and hence the
restriction of~($\dag$) to $\wt{D}_{i}$ and $E_j$ is exact. Thus, since
$f^*\cF|_{E_j}$ is trivial,
$$\,\,{1\over r(\cK)}c_1(\cK)\cdot E_j\le {1\over r(\cF)}c_1(f^*\cF)\cdot
E_j\,.$$
By~($\sharp$) and~($\flat$)  we conclude that there exists $1\le i\le s$
such that
$$\,\,{1\over r(\cK)}c_1(\cK)\cdot \wt{D}_{i}
\ge {1\over r(\cF)}c_1(f^*\cF)\cdot \wt{D}_{i}\,.$$
 Since $0<r(\cK)<r(f^*\cF)$, we conclude that $f^*\cF|_{D_{i}}$ is
not stable.
\qed

\proclaim (1.20) Proposition.
Let $C\in|nH|$ be a smooth curve. Let $X\ss\cMx$ be a closed subset such that
$F|_C$ is
locally-free and stable for all $[F]\in X$. Let $k$ be a positive integer such
that there exist smooth curves in $|knH|$ (e.g.~$k\gg 0$). Then there exists a
smooth $D_k\in|knH|$ such that $F|_{D_k}$ is stable for all $[F]\in X$.

\pf
Since $C\in |nH|$ it follows from our hypothesis that $F$ is $\mu$-stable
for all $[F]\in X$. Thus    we can cover $X$ by   subsets $X_i$, open in the
analytic topology,  so that there exists a tautological sheaf on each $S\tm
X_i$. For convenience  of exposition we will assume that  these local
tautological sheaves fit together to give a tautological sheaf on $S\tm X$;
however, as the reader will readily check, the proof works in general.   Now
let $U_k\ss|knH|$ be the open dense  subset parametrizing smooth
curves. Let
$$\,\,\wt{Z}_k:=\{([D],[F])\in U_k\tm X|\hb{ $F|_D$ is not stable}\}\,.$$
By openness of stability, $\wt{Z}_k$ is closed in $U_k\tm X$. Let $Z_k$ be
the image of $\wt{Z}_k$ under the projection $U_k\tm X\to U_k$. Since $X$
is closed in $\cMx$, it is proper, and hence $Z_k$ is closed in $U_k$. We must
show that $Z_k\not=U_k$. (Of course, if $k=1$ this is true by hypothesis.) The
proof is by contradiction, so we assume $Z_k=U_k$. Let
$$R_1,\ldots,R_k\in U_1-Z_1$$
be $k$ distinct curves. Set
$$\,\,C_0:=R_1+\cdots +R_k\,.$$
Since $|knH|\tm X\to |knH|$ is proper, and since $Z_k=U_k$, there exists a
smooth connected curve $B$, a point $0\in B$ and a map
$$\,\,g=(g_1,g_2)\cl  B\to |knH|\tm X\,,$$
with the following properties:
\msk
\item{1.} $g_1(B^0)\ss U_k$, where $B^0:=(B-0)$.
\item{2.} $g_1^{-1}([C_0])=0$.
\item{3.} The family  $\pi\cl\cC\to B$, obtained pulling back  by $g_1$
the tautological family of curves parametrized by $|knH|$, satisfies
Items~(1) and~(2) preceding Lemma~(1.19).
\item{4.} Let $\psi\cl\cC\to S\tm X$ be given by $(\psi_1,\psi_2)$, where
$\psi_1\cl
\cC\to S$ is the natural map, and $\psi_2:=g_2\circ\pi$. Let
$\cF:=\psi^*\cE$,where
$\cE$ is a tautological family on $S\tm X$. Then the
restriction $\cF|_{\pi^{-1}(b)}$ is non-stable for all $b\not=0$.
\msk
\n
Thus we can apply Lemma~(1.19) to the bundle $\cF$ on $\cC$. We
conclude that  there exist $1\le i\le k$ and $[F]\in X$ such that $F|_{R_i}$
is not stable. This is absurd because $[R_i]\in(U_1-Z_1)$.
\qed

\msk
\n
Assume that $C$, $X$ are as in the hypotheses of Proposition~(1.20). Let
$D_k\in|knH|$ be a smooth curve as in the  proposition. Let
$$\rho_k\cl X\to \cM(D_k;\xi)$$
be the morphism given by restriction. Let $\T_k$ be the theta-divisor  on
$\cM(D_k;\xi)$.

\proclaim (1.21) Lemma.
Let notation be as above. Then there exists a positive rational $\l_k$ such
that
$$\,\,c_1(\rho_k^*\T_k)=\l_kc_1(\rho_1^*\T_1)\,.$$

\pf
First let's assume that there exists a tautological sheaf $\cE$ on $S\tm X$.
Let $\cE_k$ be
its restriction to $D_k\tm X$; by our hypotheses it is a vector bundle.  Let
$M_k$ be a vector
bundle on $D_k$ such that
$$\chi(\cE|_{D_k\tm\{x\}}\ot M_k)=0\eqno(*)$$
for $x\in X$. (As is easily verified, such an $M_k$ always exists.)  Then
Grothendieck-Riemann-Roch (see~[O] for the computation in the rank-two case)
gives
$$\,\,c\left(\pi_{X,!}\left(\cE_k\ot\pi_C^*M_k\right)\right)
=\pi_{X,*}\left[ch\left(\cE_k\ot\pi_C^*M_k\right)TdC\right]\,.$$
Considering the degree-one components of both sides of the above equality, and
using~($*$),
one gets
$$\,\,-c_1\left(\det\pi_{X,!}\left(\cE_k\ot\pi_C^*M_k\right)\right)
=r(M_k)\cdot\pi_{X,*}
\left(c_2(\cE_k)-{r(\cE_k)-1\over 2r(\cE_k)}c_1(\cE_k)^2\right)\,.$$
Now assume that the  rank of $M_k$  is minimal. Then the left-hand side of the
above equality is identified with the first Chern class of
$\rho_k^*\T_k$ (see~[DN]), while
the right-hand side equals the slant product
$$\,\,r(M_k)\cdot
\left(c_2(\cE)-{r(\cE)-1\over 2r(\cE)}c_1(\cE)^2\right)/[knH]\,.$$
This proves the lemma (with $\l_k=k\cdot rk(M_k)$) under the assumption that
there is a
tautological sheaf on $S\tm X$. In general, by Theorem~(A.5) in~[Mu] there
exists a {\it quasi-tautological} sheaf $\cF$ on $S\tm X$, i.e.~such  that for
$[F]\in X$, the restriction $\cF|_{S\tm\{[F]\}}$ is isomorphic to $F^{\op\s}$
for some positive integer $\s$. Then one can repeat the  proof above with
$\cE$ replaced by $\cF$. \qed

\msk
\n
{\it Proof of Proposition~(1.18).}
\hskip 2mm
By Serre's vanishing Theorem, if $k\gg 0$ then for all
$[F_1],[F_2]\in X$ we have
$$\,\,H^1\left(F_1^*\ot F_2(-knH)\right)=0\,.\eqno(\star)$$
By Proposition~(1.20) there exists a smooth $D_k\in |knH|$ such that
$F|_{D_k}$ is stable for all $[F]\in X$. By~($\dag$) the  restriction map
$\rho_k$ is an injection. Since $\T_k$ is ample~[DN], we conclude that
$$\,\,\left(\rho_k^*\T_k\right)^{\dim X}>0\,.$$
Proposition~(1.18) follows at once from the above inequality and Lemma~(1.21).
\bsk

\n
{\bf Proof of Propositions~(1.1)-(1.2).}
\msk
\n
The proof is by contradiction, so we assume that $F$ is locally-free for all
$[F]\in X$.  Let $C\in |nH|$ be a smooth curve.    We
start by showing that there exists $[E]\in X$ such that $E|_C$  is not stable.
This again we  prove by
 reductio ad absurdum. Clearly we can assume $X$ is irreducible. If $E|_C$ is
stable for all $[E]\in X$, then we have the restriction morphism
$$\,\,\rho\cl X\to\cM(C;\xi)\,.$$
Let $g$ be the genus of $C$. By adjunction
$$\,\,g-1={1\over 2}H^2n^2+{1\over 2}K\cdot Hn\,,\eqno(\star)$$
and hence
$$\,\,\dim\cM(C;\xi)=\left(\rx^2-1\right)(g-1)
={1\over 2}\left(\rx^2-1\right)\left(H^2n^2+K\cdot Hn\right)\,.$$
By  Item~(1) of Proposition~(1.1) we see that  $\dim
X>\dim\cM(C;\xi)$, and thus
$$\,\,\left(\rho^*\T\right)^{\dim X}=0\,.$$
This contradicts Proposition~(1.18). Hence we conclude that there exists
$[E]\in X$ such that
$E|_C$ is not stable. Let $X_C\ss X$ be the (closed) subset parametrizing
sheaves whose restriction to $C$ is  not-stable; we have just proved that
$X_C\not=\es$. Let $X_C^{\mu}\ss X_C$ be the subset parametrizing
$\mu$-stable sheaves. We claim that  $X_C^{\mu}\not=\es$.   Suppose that
the ``original'' $E$ (with $E|_C$ not stable) is not $\mu$-stable. Let $\cE$ be
the family of sheaves on $S$ parametrized by $Def^0(E)$. By Luna's \'etale
slice Theorem the map
$$\l\cl Def^0(Gr E)\to \cMx$$
induced by $\cE$ is surjective onto a neighborhood of $[E]$. Hence
$$\,\,\dim \left(\l^{-1}X\right)\ge \dim X\,.\eqno(*)$$
Obviously
$$\,\,\l^{-1}X_C=\left\{x\in \l^{-1}X\ |\ \hb{$\cE_x|_C$ is not
stable}\right\}\,.$$
By Proposition~(5.47) we have
$$\,\,\dim \left(\l^{-1}X_C\right)\ge \dim\l^{-1}X-
{\rx^2\over 4}g\,.$$
By~($\star$),  Inequality~($*$) and Item~(2) of Proposition~(1.1)  we
conclude that
$$\,\,\dim\left(\l^{-1}X_C\right) >(2\rx-1)\Dx+\e(\rx,S,H)\,.$$
Hence by Proposition~(5.40)  there exists $x\in\l^{-1}X_C$ such that
$\cE_x$ is $\mu$-stable. Then $\l(x)\in X_C^{\mu}$, and hence
$X_C^{\mu}\not=\es$.  Let $[E']\in X_C^{\mu}$;
since $E'$ is stable  there is a neighborhood (in the analytic topology)  of
$[E']$ in $\cMx$   parametrizing a  tautological family $\cF$. Applying again
Proposition~(5.47), this time  to $\cF$,  we get
$$\,\,\dim X_C^{\mu}\ge \dim X-{\rx^2\over 8}(H^2n^2+K\cdot Hn)
-{\rx^2\over 4}\,.\eqno(\dag)$$
This implies that, under the hypotheses of~(1.1), there exists $[F]\in
X_C^{\mu}$ such that $[F]\notin\cMx(\a_C)$ (by Item~(3)) and, that under the
hypotheses of~(1.2), there exists $[F]\in
X_C^{\mu}$ such that $[F]\notin\cMx^C(\a_C)$ and
$[F]\notin\cMx\left(K\cdot H/2\sqrt{H^2}\right)$ (by~(1.3) and~(1.4)).
(Here $\a_C$ is as in Proposition~(1.6).) Hence in the latter case
Corollary~(5.8) gives
$$\,\,h^0(F,F\ot K)^0\le h^0(2K)\,.$$
The above properties of $F$ together with Item~(4) of~(1.1) or
Inequality~(1.5) of~(1.2) show that the hypotheses of~(1.6) are
satisfied by $F$ and $X$. By Proposition~(1.6) we conclude that $\partial
X\not=\es$. This  contradicts the assumption $\partial X=\es$;  thus we
have  proved  Propositions~(1.1)-(1.2).
\bsk
\bsk

\n
{\xlbf 2. Closed subsets not intersecting the boundary.}
\msk
\n
In this section we will prove Theorem~A,  by applying
Propositions~(1.1) and~(1.2). As a further application of this last
proposition we will give explicit conditions ensuring that every irreducible
component of $\cMx$ has non-empty boundary (Proposition~(2.11)), when
$\rx=2$. This last result will be the key ingredient in
the proof of Theorem~C. To simplify notation we will set
$\e=\e(\rx,S,H)$, $e_K=e_K(\rx,S,H)$, and $\chi=\chi(\cO_S)$.
\bsk

\n
{\bf Proof of Theorem~A for $\rx>2$.}
\msk
\n
Let
$$\eqalign{\psi_1(\xi,n):= & {1\over 2}(\rx^2-1)H^2n^2+{1\over
2}(\rx^2-1)K\cdot Hn\,,\cr
\psi_2(\xi,n):= & {1\over 8}\rx^2H^2n^2+{1\over 8}\rx^2K\cdot Hn+
{1\over 4}\rx^2+(2\rx-1)\Dx+\e\,,\cr
\psi_3(\xi,n):=& (2\rx-1)\Dx+(2\rx^3-{39\over 8}\rx^2+4\rx-1)H^2n^2
+\left[{\rx^3\over 2}|K\cdot H|+\rx^5H^2\right]n+\cr
& +{\rx^4\over 5}{(K\cdot H)^2\over H^2}+{\rx^7\over 2}H^2
{\rx^2\over 2}{(K\cdot H+\rx^2H^2)^2\over H^2}+\rx^2|\chi|+{\rx^3\over
8}|K^2|\,,\cr
\psi_4(\xi,n):=& 2\rx\Dx-{1\over 2}(r_{\xi}-1)H^2n^2+{1\over
2}(r_{\xi}-1)K\cdot Hn
-(r_{\xi}^2-1)\chi+e_K+{r_{\xi}^2\over 4}\,.}$$
If $i=1,4$ then $\psi_i(\xi,n)$ equals the right-hand side of the inequality
in Item~(i) of Proposition~(1.1). Clearly $\psi_2(\xi,n)$ bounds from
above the right-hand side of Item~(2) of~(1.1). Finally
Proposition~(5.9) and easy estimates show that $\psi_3(\xi,n)$ is an
upper bound for the right-hand side of Item~(3) in the same proposition.
Thus if
$$\,\,\dim X>\max\{\psi_1(\xi,n),\psi_2(\xi,n),
\psi_3(\xi,n),\psi_4(\xi,n)\}\,,$$
for some integer $n\ge 1$, then we conclude that $\partial X\not=\es$.   As is
easily
checked
$$\max\{\psi_1(\xi,n),\psi_2(\xi,n)\}\le\psi_3(\xi,n)\eqno(2.1)$$
 for $n\ge 1$. Hence we have

\proclaim (2.2).
Keep notation as above. If
$$\dim X>\max\{\psi_3(\xi,n),\psi_4(\xi,n)\}$$
for some integer $n\ge 1$, then $\partial X\not=\es$.

\n
If $\Dx$ is sufficiently large, then the minimum of
$\max\{\psi_3(\xi,n),\psi_4(\xi,n)\}$ for positive $n$
is achieved by the solution of
$$\,\,\psi_3(\xi,n)=\psi_4(\xi,n)\,.\eqno(\dag)$$
So let $x_0$ be the positive root of the equation in $n$
$$\,\,(2\rx-1)\Dx+(2\rx^3-{39\over 8}\rx^2+4\rx-1)H^2n^2=
2\rx\Dx-{\rx-1\over 2}H^2n^2\,,$$
obtained replacing the two sides of~($\dag$) by their dominant terms (that is
dominant for
$\Dx$ and $n$ large). Thus
$$\,\,x_0=\sqrt{{\rho(\rx)\over H^2}}\sqrt{\Dx}\,.$$
Set $n_0:=[x_0]$. We will prove Theorem~A by applying~(2.2)
with $n=n_0$. By the discussion above  this choice of $n$ is almost
optimal if $\Dx$ is large (and with this choice the computations are
relatively simple).

\proclaim Lemma.
Let $X\ss\cMx$ be a closed irreducible subset such that
$$\,\,\dim X>\max\{\psi_3(\xi,x_0),\psi_4(\xi,x_0-1)\}\,.\eqno(2.3)$$
Then $\partial X\not=\es$.

\pf
First notice that, since $\Dx\ge\D_0$, we have $x_0\ge 1$, and hence $n_0$ is a
positive
integer. If $\xi$ is fixed, the function $\psi_3(\xi,n)$ is increasing for
positive $n$.
Thus by~(2.3) we have
$$\,\,\dim X>\psi_3(\xi,n_0)\,.\eqno(*)$$
Now let's show that
$$\,\,\dim X>\psi_4(\xi,n_0)\,.\eqno(\dag)$$
First we will prove that $\dim X>\psi_4(\xi,x_0-1)$ implies that
$$\,\,x_0-1\ge  {K\cdot H\over2H^2}\,,\eqno(\star)$$
or, in other words, if the hypotheses of Theorem A are  satisfied by some
$X\ss\cMx$,
then
$$\,\,\Dx\ge\rho^{-1}H^2\left(1+{K\cdot H\over 2H^2}\right)^2\,.$$
For this observe that, if $\xi$ is fixed, then the unique critical point
of the concave-down quadratic polynomial $\psi_4(\xi,n)$ is given by $(K\cdot
H)/2H^2$.
Hence if~($\star$) does not hold then, since $0\le(x_0-1)$, we have
$$\,\,\psi_4(\xi,x_0-1)\ge\psi_4\left(\xi,{K\cdot H\over H^2}\right)>
2\rx\Dx-(r_{\xi}^2-1)\chi+e_K\,.$$
By Inequality~(0.2) we conclude that all points of $X$
parametrize non-stable sheaves. On the other hand, by~($*$) and~(2.1)
we have $\dim X>\psi_2(\xi,n_0)$. As is easily checked this implies that
$X$ satisfies the hypotheses of Corollary~(5.45); thus by this same
corollary we get a contradiction. We conclude that~($\star$) holds.
Now~($\dag$) follows at once from~(2.3), ($\star$) and the fact that
$\psi_4(\xi,n)$ is decreasing for $n\ge(K\cdot H)/2H^2$.
\qed

\msk
\n
Now we can finish the proof of Theorem~A for  $\rx>2$. A
straightforward computation gives
$$\eqalign{\psi_3(\xi,x_0)=& \l_2\Dx+
\sqrt{\rho}\left[{\rx^3\over 2}{|K\cdot H|\over
\sqrt{H^2}}+\rx^5\sqrt{H^2}\right]\sqrt{\Dx}\cr
& +{\rx^4\over 5}{(K\cdot H)^2\over H^2}+{\rx^7\over 2}H^2
+{\rx^2\over 2}{(K\cdot H+\rx^2H^2)^2\over H^2}+\rx^2|\chi|+{\rx^3\over
8}|K^2|\,,\cr
\psi_4(\xi,x_0-1)= & \l_2\Dx+\sqrt{\rho}\left[{\rx-1\over 2}{K\cdot
H\over\sqrt{H^2}}
+(\rx-1)\sqrt{H^2}\right]\sqrt{\Dx}\cr
& -{\rx-1\over 2}H^2-{\rx-1\over 2}K\cdot H-(\rx^2-1)\chi+e_K+{\rx^2\over
4}\,.}$$
As is easily checked, if $\dim X$ satisfies~(0.5) then it is greater than both
these
quantities.  By the previous lemma we conclude that $\partial X\not=\es$.
\bsk

\n
{\bf Proof of Theorem A when $\rx=2$.}
\msk
\n
The proof will be similar to the one given above, with the  difference
that instead of Proposition~(1.1) we will be using~(1.2). We set
$P_2:=h^0(2K)$.

\proclaim (2.4) Lemma.
Assume that $H$ is effective. Let $\xi$ be a set of sheaf data, with $\rx=2$.
Let $C\in |nH|$ be a smooth curve. Then
$$\eqalign{\dim\cMx^C(n\sqrt{H^2})\le & 3\Dx+2H^2n^2+(K\cdot H+2H^2+2)n\cr
& +{3(K\cdot H+H^2+1)^2\over 2H^2}+{(K\cdot H)^2\over 4H^2}
-{K^2\over 4}+3-3\chi-q_S\,.}$$

\pf
By Propositions~(5.10) and~(5.11) we have
$$\eqalign{\dim\cMx^C\left(n\sqrt{H^2}\right)\le &
\max\left\{3\a_0^2+\left((K\cdot H+2H^2+2)(H^2)^{-1/2}-n\sqrt{H^2}\right)\a_0
\right\}_{0\le \a_0\le \a}+\cr
& {3(K\cdot H+H^2+1)^2\over 2H^2}+{(K\cdot H)^2\over 4H^2}
-{K^2\over 4}+3-3\chi(\cO_S)-q_S\,.}$$
The above expression is a concave-up function of
$\a_0$, and hence its maximum is achieved at one of its end-points. A
computation shows
that its maximum is achieved for $\a_0=\a$, and equals the right-hand side of
the
inequality in the lemma.
\qed

\msk
\n
Set
$$\eqalign{\phi_1(\xi,n):= &{3\over 2}H^2n^2+{3\over 2}K\cdot Hn\,,\cr
\phi_2(\xi,n):=&{1\over 2}H^2n^2+{1\over 2}K\cdot Hn+3\Dx+\e+{3\over 2}\,,\cr
\phi_3(\xi,n):= &3\Dx+{5\over 2}H^2n^2+\left({3\over 2}K\cdot
H+2H^2+2\right)n\cr
& +{3(K\cdot H+H^2+1)^2\over 2H^2}+{(K\cdot H)^2\over 4H^2}
-{K^2\over 4}+4-3\chi-q_S\,,\cr
\phi_4(\xi,n):= &3\Dx+{1\over 2}H^2n^2+{1\over 2}K\cdot Hn+1+
\tau(S,H)\,,\cr
\phi_5(\xi,n):=& 4\Dx-{1\over 2}H^2n^2+{1\over 2}K\cdot Hn
+P_2-3\chi+1\,,}$$
where
$$\eqalign{\tau(S,H):= & {5(K\cdot H)^2\over 4H^2}+K\cdot H+
{K\cdot H\over H^2} +{3(K\cdot H+H^2+1)^2\over 2H^2}+\cr
& {(K\cdot H)^2\over 4H^2}-{K^2\over 4}+3-3\chi-q_S\quad
\hb{if $K\cdot H\ge 0$}\cr
\tau(S,H):= & 0\quad \hb{if $K\cdot H<0$.}}$$
For $i=1,2,5$ the value of $\phi_i(\xi,n)$ equals the right-hand side of the
inequality in Items~(1)-(2) of Proposition~(1.1) (with $\rx=2$), and of
Inequality~(1.5), respectively. By Lemma~(2.4) the right-hand side
of~(1.3) is  bounded above by  $\phi_3(\xi,n)$, and by  Proposition~(5.10)
$\phi_4(\xi,n)$ is an upper bound for the right-hand side of~(1.4). Thus by
Proposition~(1.2) it suffices to show that, for some  integer $n\ge 1$,
we have
$$\,\,\dim X>\max\{\phi_1(\xi,n),\phi_2(\xi,n),
\phi_3(\xi,n),\phi_4(\xi,n),\phi_5(\xi,n)\}\,.$$
As is easily checked,
$$\,\,\max\{\phi_1(\xi,n),\phi_2(\xi,n)\}\le\phi_3(\xi,n)
\quad \hb{for $n\ge 0$}\,,$$
and thus we have
\bsk

\proclaim (2.5).
Keep notation as above. If
$$\dim
X>\max\left\{\phi_3(\xi,n),\phi_4(\xi,n),\phi_5(\xi,n)\right\}\eqno(2.6)$$
for some integer $n\ge 1$, then $\partial X\not=\es$.

\n
Proceeding as in the previous case, we let $w_0$ be the positive root of
the equation in $n$
$$\,\,{5\over 2}H^2n^2+3\Dx=-{1\over 2}H^2n^2+4\Dx\,,$$
obtained by equating the two dominant terms of $\phi_3(\xi,n)$ and
$\phi_5(\xi,n)$. Explicitely
$$\,\,w_0={\sqrt{\Dx}\over \sqrt{3H^2}}\,.$$

\proclaim (2.7) Lemma.
Keep notation as above. If
$$\,\,\dim X>\max
\left\{\phi_3(\xi,w_0),\phi_4(\xi,w_0),
\phi_5(\xi,w_0-1)\right\}\,,\eqno(2.8)$$
then $\partial X\not=\es$.

\pf
 Let $n_0:=[w_0]$. Since $\Dx>\D_0(2,S,H)$, we have $n_0\ge 1$. We claim
that~(2.6) holds with $n=n_0$; by~(2.5) this will imply the
lemma. If $\xi$ is fixed, then
\msk
\item{1.} $\phi_3(\xi,n)$ is increasing for $n\ge 1$,
\item{2.} $\phi_4(\xi,n)$ is increasing for $n\ge -(K\cdot H)/2H^2$, and
\item{3.} $\phi_5(\xi,n)$ is decreasing for $n\ge (K\cdot H)/2H^2$.
\msk
\n
Since $n_0\ge 1$, $\phi_3(\xi,n)$  is
increasing for $n\ge n_0$, and since $\Dx>\D_0$,
$\phi_5(\xi,n)$ is decreasing for $n\ge(w_0-1)$.   Hence~(2.8) implies  that
$$\,\,\dim X>\max\left\{ \phi_3(\xi,n_0),\phi_5(\xi,n_0)\right\}\,.$$
If $K\cdot H\ge 0$ then, by Item~(2) above, Inequality~(2.8) also implies
$\dim X>\phi_4(\xi,n_0)$, so we are done. If $K\cdot H<0$, then as is easily
checked $\phi_4(\xi,n)\le\phi_3(\xi,n)$ for all $n\ge 1$
(use~(5.5) to take care of the contribution of $\chi$  in $\phi_3$).
Hence also in this case $\dim X>\phi_4(\xi,n_0)$.
\qed

\msk
\n
Now one finishes the proof of Theorem~A in the rank-two case by checking
that if~(0.5) holds, then  the hypotheses of the above lemma are
satisfied.
\bsk

\n
{\bf Another application of Proposition~(1.2).}
\msk
\n
Our goal in this subsection is to determine an effective $\D_1$ with the
property that, if
$\Dx>\D_1$, then every closed subset of $\cMx$ whose dimension is at least the
expected
dimension of $\cMx$ has non-empty boundary. One  such  lower bound can be
obtained by
applying Theorem~A.  However, while Theorem~A  provides the best
``asymptotic'' result of its kind obtainable from Proposition~(1.1), it is
not sharp for $\Dx$ small; hence we  proceed differently.    We will
limit ourselves to  rank-two sheaves.  Let $z_0$ be the  positive root of
the equation in $n$
$$\,\,4\Dx-3\chi=\phi_5(\xi,n)\,,\eqno(2.9)$$
i.e.
$$\,\,z_0={1\over 2H^2}
\left[K\cdot H+\sqrt{(K\cdot H)^2+8H^2(P_2+1)}\right]\,.\eqno(2.10)$$

\proclaim (2.11) Proposition.
Assume that $H$ is effective. Keeping notation as above, let
$$\eqalignno{\D_2(S,H):= &(4K\cdot H+7H^2+2)z_0+6H^2+{9\over 2}K\cdot H
+3{K\cdot H\over H^2} +{7\over 4}{(K\cdot H)^2\over H^2}+{3\over 2H^2}+\cr
& 14+5P_2(S)-{K^2\over 4}-q_S\,, & (2.12)\cr
n_0:= & \hb{\rm the least positive integer such that $n_0>z_0$.} &(2.13)}$$
Let   $\xi$ be a set of sheaf data, with   $\rx=2$. Assume that
$\Dx>\D_2(S,H)$ and that the linear system $|n_0H|$ contains  a smooth
curve. Then the following hold:
\msk
\item{1.} If $Y$ is an irreducible component of   $\cMx$, then the
generic point of $Y$ parametrizes a $\mu$-stable sheaf, and hence
$$\,\,\dim Y\ge 4\Dx-3\chi(\cO_S)\,.$$
\item{2.} If $X\ss\cMx$ is a closed irreducible subset such that
$$\,\,\dim X\ge 4\Dx-3\chi(\cO_S)\,,\eqno(2.14)$$
 then $\partial X\not=\es$.

\pf
We start by proving Item~(2).
We will show that
$$\,\,4\Dx-3\chi>\max
\left\{\phi_3(\xi,n_0),\phi_4(\xi,n_0),\phi_5(\xi,n_0)\right\}\,.\eqno(*)$$
Item~(2) will then follow from~(2.14) and~(2.5).
Since $z_0$ is the positive root of~(2.9), and since $\phi_5(\xi,n)$ is a
concave-down function of $n$, we have
$$\,\,4\Dx-3\chi>\phi_5(\xi,n_0)\,,$$
for any $\Dx$ (i.e.~even if $\Dx\le\D_2$).
Now let's first examine the case $K\cdot H\ge 0$. In this case both
$\phi_3(\xi,n)$ and $\phi_4(\xi,n)$ are increasing functions for
$n\ge 1$ (see the proof of Lemma~(2.7)). Thus it suffices to show that
$$\eqalignno{4\Dx-3\chi>&\phi_3(\xi,z_0+1)\,,&(\dag)\cr
4\Dx-3\chi>&\phi_4(\xi,z_0+1)\,.}$$
Computing, one gets that the two above inequalities are equivalent to
$$\eqalignno{\Dx> & \D_2(S,H) & (\star)\cr
\Dx> & (K\cdot H+H^2)z_0+2H^2+{9\over 2}K\cdot H+
4{K\cdot H\over H^2}+3{(K\cdot H)^2\over H^2}+{3\over 2H^2}+
 8+P_2-{K^2\over 4}-q_S\,,}$$
respectively. As is easily checked using~(2.10), the first inequality implies
the second.
This proves~($*$) if $K\cdot H\ge 0$. Now let's assume $K\cdot H<0$.
In this case  $\phi_3(\xi,n)$ is again  increasing for $n\ge
1$ and $\phi_4(\xi,n_0)\le\phi_3(\xi,n_0)$ (see the proof of
Lemma~(2.7)). Hence  it suffices to show that~(\dag) holds. Since
this inequality is equivalent to~($\star$) we are done.
Now let's prove Item~(1). As is easily checked one has
$$\,\,\phi_3(\xi,n_0)>3\Dx+\e\,,$$
and hence
$$\,\,4\Dx-3\chi>3\Dx+\e\,.$$
Item~(1) follows from the above inequality and Corollary~(5.46).
\qed

\msk
\n
We wish to rewrite $\D_2(S,H)$ when
$$c_1(K)=kc_1(H)$$
for some rational positive $k$. Let
$$\eqalign{N_2(S,H):= & {17\over 2}+6\sqrt{1+{8(\chi+1)\over 9K^2}}+
\left[8+{4\over H^2}+\left({21\over 2}+{3\over H^2}\right)
\sqrt{1+{8(\chi+1)\over 9K^2}}\right]k^{-1}+\cr
& \left(6+{14\over H^2}+{3\over 2(H^2)^2}\right)k^{-2}\,.}\eqno(2.15)$$
A straightforward computation  gives

\proclaim (2.16).
Let $(S,H)$ be a polarized surface. Assume that $c_1(K)=kc_1(H)$ for
a rational positive $k$. (In particular $K$ is ample.) Then
$$\,\,\D_2(S,H):=N_2(S,H)K^2+5\chi-q_S\,.$$

\bsk

\n
{\bf Comments.}
\msk
\n
Theorem~A naturally raises the question: what is the maximum
dimension of closed subsets of $\cMx$ which do not
intersect $\partial\cMx$ ? The following proposition gives a lower bound
for this maximum which can be contrasted with the upper bound provided
by Theorem~A.

\proclaim  Proposition.
Let $S$ be a surface such that $H^{1,1}_{\bf Z}(S)=\ZZ c_1(C)$, where $C$ is
a curve of (arithmetic) genus $g$. Fix an integer $r\ge 2$ and, for
$d\in\ZZ$, let
$$\,\,\xi(d)=(r,-[C],d)\,.$$
If  $d$ is sufficiently large there exists a closed subvariety
$X_d\ss\cM_{\xi(d)}$ such that
$$\eqalignno{\partial X_d= & \es\,,\cr
\dim X_d\ge & (r-2)d-(r-2)(r+g-1)-\dim Aut(C)\,. & (2.17)}$$

\pf
Let $L$ be a non-special line bundle on $C$.
Set $d=\deg(L)$, $n:=h^0(L)-1$. We assume that:
\msk
\item{1.} The complete linear system $|L|$ defines an  embedding
$$\,\,C\hra \PP^n=\PP\left(H^0(L)^*\right)\,.$$
\item{2.} $n+1\ge r$ and $d> C\cdot C$.
\msk
\n
 Given an $r$-dimensional subspace
$V\ss H^0(L)$ with no base-points, let $F_V$ be the sheaf on $S$ fitting
into the exact sequence
$$\,\,0\to F_V\to V\ot\cO_S\brel e_V\over\lra
\i_*\cO_C(L)\to 0\,,\eqno(*)$$
where $e_V$ is the evaluation map, and $\i\cl C\hra S$ is the inclusion.
Clearly $F_V$ is locally-free.  The
Chern classes of $F_V$ are given by
$$\eqalign{c_1(F_V)= & -[C]\,,\cr
c_2(F)= & d\,.}\eqno(\dag)$$
As is easily checked, it follows from
$H^{1,1}_{\bf Z}(S)=\ZZ c_1(C)$  that $F_V$ is $\mu$-stable.
One can  identify  the set of base-point free $r$-dimensional subspaces of
$H^0(L)$ with  the set $U_C$ of $(n-r)$-dimensional linear
subspaces of $\PP^n$ not intersecting $C$ (embedded
by $|L|$). Let $\PP^{n-2}\ss\PP^n$ be disjoint from $C$, and set
$$\,\,B_d:={\bf Gr}\left(n-r,\PP^{n-2}\right)\,.$$
Then  $B_d$ is a projective subset of $U_C$, and
$$\,\,\dim B_d=(r-2)(d-g-r+1)\,.\eqno(\star)$$
Clearly one can construct  a family $\cF$ of vector bundles on
$S$ paramatrized by $B_d$, with the property that if $[V]\in B_d$, then
$$\,\,\cF|_{S\tm[V]}\cong F_V\,.$$
Since the Chern classes of $F_V$ are given by~($\dag$), and since $F_V$ is
stable for all $V$, the family $\cF$ defines a morphism
$$\,\,\vf\cl B_d\to \cM_{\xi(d)}\,.$$
Let $X_d:=\vf(B_d)$. Clearly $X_d$ is closed and $\partial X_d=\es$. Now
let's check that~(2.17) holds. Dualizing~($*$) one gets
$$\,\,0\to V^{*}\ot\cO_S\to F_V^{*}\to\cO_C(C)\ot L^{*}\to 0\,.$$
By Item~(2) above we have $V^{*}\cong H^0(F_V^{*})$, and hence the
isomorphism class of $F_V$ determines $V$ up to isomorphism.
Formula~(2.17) follows at once from this and Equation~($\star$).
\qed
\bsk
\bsk

\n
{\xlbf 3. Moduli of bundles with twisted endomorphisms.}
\msk
\n
In this section we will prove Theorems~B, C, and their corollaries. First
some notation. If $X\ss\cMx$, we let $\ov{X}\supset X$ be the closure of $X$
in $\cMx$,   and  $X^{\mu}\ss X$ be the subset parametrizing $\mu$-stable
sheaves.   To
simplify notation we will set
$$\,\,\pal\ov{X}=\pal\left(\ov{X}\right)\quad\pal X^{\mu}:=(\pal X)^{\mu}
\quad\pal\ov{X}^{\mu}:=\left(\pal\ov{X}\right)^{\mu}\,.$$
\bsk

\n
{\bf The double-dual construction.}
\msk
\n
The basic idea in the proof of Theorems~B and~C
is to consider the double-duals of sheaves parametrized by $\partial
\ov{W}_{\xi}^L$, or by a non-good irreducible component of $\cMx$
respectively. In this subsection we will establish some results concerning
this ``double-dual construction''.  Assume that $X\ss\cMx$ and that $\pal
X^{\mu}\not=\es$.  Let $[F]\in\pal X^{\mu}$;  we have the canonical exact
sequence
$$\,\,0\to F\to F^{**}\brel\psi_F\over\lra Q_F\to 0\,,\eqno(3.1)$$
where $Q_F$ is an Artinian sheaf of finite  length
$$\,\,\ell\left(Q_F\right)=h^0\left(Q_F\right)>0\,.$$
Since $F$ is $\mu$-stable, so is $F^{**}$, and hence it determines a point
$[F^{**}]\in\cM_{\xi'}$, where
$$\,\,\xi'=(\rx,\det_{\xi},c_2(\xi)-\ell(Q_F))\,.\eqno(3.2)$$
The sheaves $F^{**}$, for $[F]$ varying in $\pal X^{\mu}$, do not fit together
to give a family of sheaves (their second Chern classes might vary). However
there is a (maximal) stratification of $\pal X^{\mu}$ by locally closed
subsets with the property that the double duals of sheaves parametrized by
points of the same stratum fit together to give (locally) a family of vector
bundles. We call this the {\it double-dual stratification} of $\pal X^{\mu}$.
We will be interested in the open strata. First we give a lower bound for their
dimension.

\proclaim (3.3) Proposition.
Let $\cF$ be a family of rank-$r$ torsion-free sheaves on $S$
parametrized by an equidimensional variety $B$. Let $\partial B\ss B$ be the
subset parametrizing singular sheaves. Then $\partial B$ is closed, and
$$\,\,\cod\left(\partial B,B\right)\le r-1\,.$$

\pf
Since $\cF$ is a family of rank-$r$ torsion-free sheaves on a smooth
surface it has a short locally-free resolution
$$\,\,0\to \cE_1\brel\phi\over\lra \cE_0\to \cF\to 0\,.$$
Let ${\rm D}(\phi)\ss S\tm B$ be the degeneracy locus of $\phi$
(i.e.~the locus where $\phi$ drops rank). Clearly ${\rm D}(\phi)$ is
closed, and since $r(\cE_0)-r(\cE_1)=r$ we have
$$\,\,\cod\left({\rm D}(\phi),S\tm B\right)\le r+1\,.$$
The result follows because $\pal B=\pi_B\left({\rm D}(\phi)\right)$.
\qed

\proclaim (3.4) Corollary.
Let $X\ss\cMx$ be a locally closed equidimensional subset. Assume that
$\pal X^{\mu}\not=\es$. Let $Y\ss \pal X^{\mu}$ be an irreducible
component of an open stratum of the double-dual stratification. Then
$$\,\,\dim Y\ge \dim X-(\rx-1)\,.$$

\n
Let $Y$ be as in the above corollary, and set
$$\,\,Y^{**}:=\left\{[F^{**}]\ |\ [F]\in Y\right\}\,.$$
Thus $Y^{**}\ss\cM_{\xi'}$, where $\xi'$ is given by~(3.2). We will
relate the dimensions of $X$ and $Y^{**}$. For this we need to consider
certain Quot-schemes. If $E$ is a vector bundle on $S$, and $\ell$ is a
positive integer, let $Quot(E;\ell)$ be the parameter space for quotients of
$E$ of finite length equal to $\ell$. Let $Quot_0(E;\ell)$ be the (open)
subset parametrizing quotients
$$\,\,E\to \bigoplus_{i=1}^{\ell}\CC_{P_i}\,,$$
where $\CC_{P_i}$ is the skyscraper sheaf at $P_i$. The following
result is due to Li~[Li] (the proof there  is given for $r(E)=2$, but in fact
it carries over to any rank, see~[GL2]):

\proclaim (3.5) Theorem (Li).
Let notation be as above. Then $Quot_0(E;\ell)$ is dense in $Quot(E;\ell)$.
In particular  $\dim Quot(E;\ell)=(r(E)+1)\ell$.

\proclaim  (3.6) Corollary.
Let $X\ss\cMx$, and assume that $\pal X^{\mu}\not=\es$. Let $Y\ss \pal
X^{\mu}$ be an irreducible component of an open stratum of the double-dual
stratification. Let $\ell:=\ell(Q_F)$, where $[F]\in Y$. Then
$$\,\,\dim Y^{**}\ge \dim X-(\rx+1)\ell-\rx+1\,.\eqno(3.7)$$
If~(3.7) is an equality,  the following holds: Let $[E]\in Y^{**}$, and let
$\phi\in Quot_0(E;\ell)$ be a $\ul{\rm generic}$ quotient. Then the sheaf
$F_{\phi}$ fitting into the exact sequence
$$0\to F_{\phi}\to E\brel \phi\over\lra
\bigoplus_{i=1}^{\ell}\CC_{P_i}\to 0$$
is parametrized by a point of $Y$.

\pf
Let $Quot(Y^{**};\ell)\to Y^{**}$ be the relative $Quot$-scheme, with fiber
$Quot(E;\ell)$ over $[E]\in Y^{**}$, and let $Quot_0(Y^{**};\ell)$ be the
open subset with fiber $Quot_0(E;\ell)$ over $[E]$. We have an injection
$$\,\,f\cl Y\hra Quot(Y^{**};\ell)\,,$$
mapping $[F]$ to the canonical quotient $\psi_F$ (see~(3.1)). By
Theorem~(3.5) we conclude that
$$\,\,\dim Y\le \dim Y^{**}+(\rx+1)\ell\,.\eqno(*)$$
Inequality~(3.7) follows from the above inequality and Corollary~(3.4).
Now suppose that~(3.7) is an equality then we must have equality also
in~($*$). By Theorem~(3.5) we conclude that $f(Y)\cap Quot_0(Y^{**};\ell)$ is
dense in  $Quot_0(Y^{**};\ell)$; this proves the second statement of the
corollary.
\qed

\msk
Now let $L$ be a line bundle on $S$. If $X\ss\cMx$ is an irreducible
locally closed subset with $X^{\mu}\not=\es$, we set
$$\,\,h_L(X):=\min\left\{h^0(F,F\ot L)^0|
\ [F]\in X^{\mu}\right\}\,.$$
By semicontinuity of cohomology dimension,   if $[F]\in X^{\mu}$ is a
generic point  then $h^0(F,F\ot L)^0=h_L(X)$.  The following
proposition containes the observation that will allow us to deduce
Theorems~B and~C from Theorem~A.

\proclaim (3.8) Proposition.
Let notation be as above. Let $X\ss\cMx$ be  a locally closed
irreducible subset such that $\pal X^{\mu}\not=\es$.   Let
$Y\ss\pal X^{\mu}$ be an irreducible component of an open stratum of  the
double-dual stratification, and let $\ell:=\ell(Q_F)$ for $[F]\in Y$. Then:
\msk
\item{1.} $h_L(Y^{**})\ge h_L(X)$.
\item{2.} $\dim Y^{**}\ge \dim X-(2\rx-1)\ell-1$.
\item{3.} If $h_L(X)>0$ then one at least of the  inequalities In
Items(1)-~(2) is strict.

\pf
Let $[E]\in Y^{**}$ be such that $h^0(E,E\ot L)^0=h_L(Y^{**})$. Let
$[F]\in Y$ be a sheaf such that $F^{**}\cong E$. There exists  a natural
injection
$$\,\,\d\cl H^0(F,F\ot L)^0\hra H^0(E,E\ot L)^0\,.\eqno(\star)$$
Item~(1) follows at once from this.   Item~(2)  follows from
Corollary~(3.6). Now let's prove Item~(3). So assume that we have
equality in Item~(2): we will prove that  the inequality of Item~(1) is strict.
Let $[E]\in Y^{**}$ be as above. Let $F_{\phi}$ be as in the statement of
Corollary~(3.6) where  $\phi$  is  an arbitrary quotient.
Letting $\d_{\phi}$ be as in~($\star$), with $F$ replaced by $F_{\phi}$, we
have
$$\,\,\d_{\phi}\left(H^0(F_{\phi},F_{\phi}\ot L)\right)=
\left\{f\in H^0(E,E\ot L)^0\ |\ f^*\phi_i=\l_i\phi_i\right\}\,,$$
where $\phi_i$ is the restriction of $\phi$ to $P_i$, and $\l_i\in\CC$.
Now suppose that $f\not=0$ and that $\phi$ is generic. Then,  since $f$  is
not a dilation,  it is not contained in the image of $\d_{\phi}$.
Hence
$$\,h^0(F_{\phi},F_{\phi}\ot L)^0<h_L(Y^{**})\quad
\hb{for $\phi$ generic.}\eqno(\dag)$$
On the other hand, if Item~(2) is an equality then so must be
Inequality~(3.7). Hence by Corollary~(3.6) $[F_{\phi}]\in Y$ for $\phi$
generic. By~($\dag$) we conclude that the inequality in Item~(1) is strict.
\qed
\bsk

\n
{\bf Proof of Theorem~B.}
\msk
\n
The proof will be by contradiction. So we assume that $\Dx>\D_1$, and that
Inequality~(0.6) is violated, for some line bundle $L$. Let  $X_0\ss
W_{\xi}^L$ be an irreducible component of maximum dimension.
By hypothesis
$$\,\,\dim X_0>\l_2\Dx+\l_1\sqrt{\Dx}+\l'_0+e_L\,.\eqno(3.9)$$
By the above inequlity  and  Theorem~A we have $\pal\ov{X}_0\not=\es$.
The following lemma will show that in fact $\pal\ov{X}_0^{\mu}\not=\es$.

\proclaim (3.10) Lemma.
Let $X\ss\cMx$ be an equidimensional locally closed subset such that:
\msk
\item{1.} $\pal X\not=\es$, and
\item{2.} $\dim X>(2\rx-1)\Dx+\e(\rx,S,H)+\rx-1$.
\msk
\n
Then $\pal X^{\mu}\not=\es$.

\pf
Let $[F]\in\pal X$. If $F$ is $\mu$-stable there is nothing to prove,
so assume $F$ is not $\mu$-stable. Let   $\cF$ be the family of sheaves on
$S$  parametrized by $Def^0(Gr F)$.  By Luna's \'etale slice Theorem, the
map
$$\l\cl Def^0(Gr F)\to\cMx$$
induced by $\cF$ is surjective onto a neighborhood of $[F]$, and hence
$$\,\,\dim\left(\l^{-1}X\right)\ge\dim X\,.\eqno(\star)$$
Let $\partial\left(\l^{-1}X\right)\ss\l^{-1}X$ be the subset parametrizing
singular
sheaves.  By Proposition~(3.3) and by Item~(2)  we conclude that
$$\,\,\dim\partial\left(\l^{-1}X\right)>(2\rx-1)\Dx+\e(\rx,S,H)\,.$$
By~(5.40)  the subset of
$\partial\left(\l^{-1}X\right)$ parametrizing
$\mu$-stable sheaves is non-empty. Since this subset is open, and since its
image is contained in $\pal X^{\mu}$, we see that $\pal X^{\mu}\not=\es$.
\qed

\msk
\n
An easy computation shows that~(3.9) together with $\Dx>\D_1$ give
$$\,\,\dim X_0>(2\rx-1)\Dx+\e(\rx,S,H)+\rx-1\,.$$
Hence by Lemma~(3.10) we have $\pal\ov{X}_0^{\mu}\not=\es$. Let
$Y_0\ss\pal\ov{X}_0^{\mu}$ be an irreducible component of an open
stratum of the double-dual stratification. Let
$$\,\,\xi_1:=(\rx,\detx,c_2(\xi)-\ell_0)\,,$$
where $\ell_0:=\ell(Q_{F_0}))$ for $[F_0]\in Y_0$. Set $X_1:=Y_0^{**}$.
Thus
$$\,\,X_1\ss\cM_{\xi_1}\,.$$
Now consider $\ov{X}_1$ and, if $\pal \ov{X}^{\mu}_1\not=\es$ continue in
the same fashion. We will get a sequence
$$\,\,X_i\ss\cM_{\xi_i}\,,\quad Y_i\ss\pal\ov{X}^{\mu}_i\,,
\quad X_{i+1}=Y_i^{**}\,,\eqno(3.11)$$
for $i=0,\ldots,n$ (with $\xi_0:=\xi$), until we reach a point when
$\pal\ov{X}_n^{\mu}=\es$.  We will show that the dimension of $X_n$ is
``too big'', and thus get a contradiction.     Let $\ell_i:=\ell(Q_{F_i})$,
where $[F_i]\in X_i$, and set
$$\,\,\ell:=\ell_0+\cdots\ell_{n-1}\,.$$

\proclaim Lemma.
Keeping notation as above, we have
$$\,\,\dim X_n\ge \dim X_0-(2\rx-1)\ell-e_L\,.$$

\pf
By Item~(1) of~(3.8) we have
$$\,\,h_L(X_0)\le h_L(X_1)\le\cdots\le h_L(X_n)\,.\eqno(*)$$
In particular, since $h_L(X_0)>0$, we have $h_L(X_i)>0$ for all $i$. Since
by~(5.7) $h_L(X_n)\le e_L$ there can be at most $(e_L-1)$ strict
inequalities in~($*$). Thus by Proposition~(3.8) again we have
$$\,\,\dim X_{i+1}\ge \dim X_i-(2\rx-1)\ell_i-\d_i\,,$$
where $\d_i=0$ or $\d_i=1$, and $\d_i$ is equal to $1$ for  at most $(e_L-1)$
values of $i$.  The result follows at once from this.
\qed

\msk
\n
The above Lemma together with~(3.9) gives
$$\eqalignno{\dim X_n > & \l_2\Dx+\l_1\sqrt{\Dx}+\l'_0-(2\rx-1)\ell\cr
= & (2\rx-1)\D_{\xi_n}+\left[\l_2-(2\rx-1)\right]\Dx+
\l_1\sqrt{\Dx}+\l'_0 & (3.12) \cr
> &
\max\left\{\l_2\D_{\xi_n}+\l_1\sqrt{\D_{\xi_n}}+\l'_0,
(2\rx-1)\D_{\xi_n}+\e+\rx-1\right\}\,. &(3.13)}$$

\proclaim Lemma.
Keeping notation as above, we have
$$\,\,\D_{\xi_n}\le\D_0\,.\eqno(3.14)$$

\pf
The proof is by contradiction, so we assume that~(3.14) is violated.
Then Inequality~(3.13) and Theorem~A give that $\pal
\ov{X}_n\not=\es$. By~(3.13) and Lemma~(3.10) we conclude that
$\pal\ov{X}_n^{\mu}\not=\es$. This contradicts the definition of $X_n$.
\qed

\msk
Now we can finish the proof of Theorem~B. We will show that
$$\,\,\dim X_n>2\rx\Dx-(\rx^2-1)\chi+e_K\,.\eqno(*)$$
{}From this one concludes as follows: Since $X_n\ss\cM_{\xi_n}$ the above
inequality together with~(0.3) implies that all sheaves parametrized by
$X_n$ are non-stable. But this is absurd by Inequality~(3.13) and by~(5.45).
Now let's prove~($*$). By~(3.12)  it suffices to show that
$$\,\,\left[\l_2-(2\rx-1)\right]\Dx+\l_1\sqrt{\Dx}+\l'_0
>\D_{\xi_n}-(\rx^2-1)\chi+e_K\,.$$
By Inequality~(3.14) it suffices to check that
$$\,\,\left[\l_2-(2\rx-1)\right]\Dx+\l_1\sqrt{\Dx}+\l'_0
>\D_0-(\rx^2-1)\chi+e_K\,.$$
This follows at once from $\Dx>\D_1$.
\qed

\bsk
\n
{\bf Proof of Corollary ${\bf B}'$.}
\msk
\n
Let $\wt{\D}_1(\rx,S,H)$ be the smallest number such that:
\msk
\item{$\bullet$} $\wt{\D}_1\ge\D_1$, and
\item{$\bullet$} if  $\Dx>\wt{\D}_1$ then
$$\eqalignno{2\rx\Dx-(\rx^2-1)\chi(\cO_S)> &
(2\rx-1)\Dx+\e(\rx,S,H)+\rx\,, &(3.15)\cr
2\rx\Dx-(\rx^2-1)\chi(\cO_S)> & \l_2\Dx+\l_1\Dx+\l'_0+e_K\,. &(3.16)}$$
\msk
\n
Set
$$\,\,\D'_1(\rx,S,H):=\wt{\D}_1(\rx,S,H)+(\rx-1)^{-1}e_K(\rx,S,H)\,.$$
As is easily checked $\D'_1$ depends only on $\rx$, $K^2$, $K\cdot H$,
$H^2$ and $\chi(\cO_S)$.

\proclaim (3.17) Claim.
Let $(S,H)$ be a polarized surface, with $H$ satisfying~(0.4). Let $\xi$ be
a set of sheaf data such that $\Dx>\wt{\D}_1(\rx,S,H)$. Then $\cMx$ is good,
and the generic point of any of its irreducible components parametrizes a
$\mu$-stable sheaf.

\pf
Let $X$ be an irreducible component of $\cMx$. Inequality~(3.15)
and Corollary~(5.46) give that the generic point of $X$ parametrizes a
$\mu$-stable sheaf. Since $\Dx>\D_1$, we conclude by Theorem~B and
Inequality~(3.16) that
$$\,\,\dim \Wx^K<2\rx\Dx-(\rx^2-1)\chi(\cO_S)\,.$$
Thus by Inequlaity~(0.2) we have
$$\,\,\dim X>\Wx^K\,,$$
i.e.~$h^0(F,F\ot K)^0=0$ for a (stable) sheaf $F$ parametrized by the generic
point of $X$. Thus $\cMx$ is good.
\qed

\msk
Now assume that $\Dx>\D'_1$. Let $X$ be an irreducible component of
$\cMx$. By~(3.17) we have $X^{\mu}\not=\es$. We will show that if $[F]\in
X^{\mu}$ is generic then $F$ is locally-free.   Assume the contrary,
i.e.~$X^{\mu}\ss\pal\cMx$.  Let $Y\ss X^{\mu}$ be an irreducible component of
an open stratum of the double-dual stratification. Then   %
$$\,\,Y^{**}\ss\cM_{\xi'}\quad c_2(\xi)-c_2(\xi')>0\,.$$
It follows from Theorem~(3.5) that
$$\,\,\dim X^{\mu}\le \dim Y^{**}+
(\rx+1)\left(c_2(\xi)-c_2(\xi')\right)\,.\eqno(*)$$
We distinguish between two cases:
\msk
\item{1.} If $\D_{\xi'}>\wt{\D}_1$, then $\dim Y^{**}\le
2\rx\D_{\xi'}-(\rx^2-1)\chi(\cO_S)$ by Claim~(3.17).
\item{2.} If $\D_{\xi'}\le\wt{\D}_1$, then $\dim Y^{**}\le
2\rx\D_{\xi'}-(\rx^2-1)\chi(\cO_S)+e_K$ by~(0.3).
\msk
\n
In both cases, Inequality~($*$) gives
$$\,\,\dim X^{\mu}<2\rx\D_{\xi}-(\rx^2-1)\chi(\cO_S)\,.$$
This is absurd because $\cMx$ is good. Hence the generic point of $X^{\mu}$
parametrizes a locally-free sheaf. This finishes the proof of
Corollary~${\rm C}'$.
\bsk

\n
{\bf Proof of Theorem C.}
\msk
\n
The essential step is provided by the following

\proclaim (3.18) Proposition.
Let $(S,H)$ be a polarized surface satisfying the hypotheses of Theorem~C.
Let $\xi$ be a set of sheaf data with $\rx=2$. Assume that
$$\,\,\Dx>\D_2(S,H)+2h^0(2K)\,.$$
If $X_0\ss\cMx$ is an irreducible component then there exists $[F]\in\pal
X_0^{\mu}$ such that
$$\,\,h^0(F^{**},F^{**}\ot K)^0=0\,.\eqno(3.19)$$

\pf
By Proposition~(2.11) we have $\pal X_0\not=\es$. As is easily checked we
have
$$\,\,4\D_2-3\chi>  3\D_2+\e+1\,.\eqno(3.20)$$
Thus by   Lemma~(3.10) we conclude that $\pal X_0^{\mu}\not=\es$. Let
$Y_0\ss\pal X_0^{\mu}$ be an irreducible component of an open stratum of
the double-dual stratification, and set $X_1:=Y_0^{**}$. Assume that
$$\,\,h_K(X_1)>0\,.\eqno(3.21)$$
We will arrive at a contradiction, and thus we will conclude that
$h_K(X_1)=0$; this will prove the proposition. Proceeding as in the
proof of Theorem~B, we construct a series of locally closed irreducible
subsets as in~(3.11); the only difference is that in  the present case we
define
$n$ by requiring that:
\msk
\item{1.} $\D_{\xi_{n-1}}>\D_2$, and
\item{2.} either $\pal\ov{X}_n^{\mu}=\es$, or $\D_{\xi_n}\le\D_2$.
\msk

\proclaim Claim.
Keeping notation as above, we have
$$\eqalignno{\D_{\xi_n}\le & \D_2\,, & (3.22)\cr
h_K(X_i)\le & h^0(2K)\quad \hb{for $0\le i\le n-1$.} & (3.23)}$$

\n
{\bf Proof of the claim.}
\hskip 2mm
Since $\Dx>\D_2$ we have $\dim X_0\ge 4\Dx-3\chi$ (Proposition~(2.11)).
Thus, by Item~(2) of Proposition~(3.8) we have
$$\,\,\dim X_i\ge 4\D_{\xi_i}-3\chi
\quad\hb{for all $0\le i\le n$.}\eqno(3.24)$$
Now assume that $\D_{\xi_n}>\D_2$. Then by the above
inequality and by Proposition~(2.11) we have $\pal\ov{X}_n\not=\es$. By
Inequality~(3.20) and Lemma~(3.10) we conclude that
$\pal\ov{X}^{\mu}_n\not=\es$. This contradicts the definition of $n$, and
thus we conclude that~(3.22) holds. In order to prove~(3.23) we need the
following easily checked inequality:
$$\eqalign{4\D_2-3\chi> & 3\D_2+3{(K\cdot H)^2\over 4H^2}
+{(K\cdot H)(K\cdot H+2H^2+2)\over 2H^2} \cr
& +{3(K\cdot H+H^2+1)^2\over 2H^2}+{(K\cdot H)^2\over 4H^2}
-{K^2\over 4}+3-3\chi(\cO_S)-q_S\,.}\eqno(3.25)$$
The right-hand side of~(3.25) equals the right-hand side of the inequality in
Proposition~(5.10) with $\a:=(K\cdot H)/2\sqrt{H^2}$, and $\Dx$ replaced
by $\D_2$. The above inequality, together with~(3.24) and Item~(1)
preceding the claim, gives
$$\dim X_i>\dim\cM_{\xi_i}\left({K\cdot H\over 2\sqrt{H^2}}\right)
\quad\hb{for $0\le i\le n-1$.}$$
By Corollary~(5.8) we conclude that~(3.23) holds.
\qed

\msk
\n
Now we will use the full strenght of Proposition~(3.8) to give a
lower bound for $\dim X_n$ which is larger than~(3.24).   By~(3.21) and
by Item~(1) of~(3.8) we have
$$\,\,0< h_K(X_1)\le h_K(X_{n-1})\le h_K(X_n)\,.$$
By~(3.23) there are at most $h^0(2K)$ strict inequalities. Thus
applying  Proposition~(3.8) repeatedly to $X_0$, $X_1$, ... one gets
$$\,\,\dim X_n\ge \dim X_0-3\left(c_2(\xi)-c_2(\xi_n)\right)-h^0(2K)\ge
3\D_{\xi_n}+\Dx-3\chi-h^0(2K)\,.$$
Since $\Dx>\D_2+2h^0(2K)$ we conclude that
$$\,\,\dim X_n> 3\D_{\xi_n}-3\chi+\D_2+h^0(2K)\,.\eqno(\star)$$
By construction the generic point of $X_n$ parametrizes a $\mu$-stable
vector bundle. Hence:
\msk
\item{1.} If $e_K(X_n)\le h^0(2K)$ then by~(0.3) we have $\dim X_n\le
4\D_{\xi_n}-3\chi+h^0(2K)$,  and
\item{2.} if $e_K(X_n)> h^0(2K)$, then by Corollary~(5.8)
$\dim X_n\le  \dim\cM_{\xi_n}\left((K\cdot H)/2\sqrt{H^2}\right)$.
\msk
\n
By~(3.22) and by~($\star$) we see that Item~(1) is
impossible. On the other hand also Item~(2) is impossible,
by~(5.10), ($\star$), (3.22) and~(3.25). Thus we have a contradiction.
This proves that~(3.21) can not hold, and hence proves the proposition.
\qed

\msk
Now we can prove Theorem~C. Let $X_0\ss\cMx$ be an irreducible
component. Let $[F]\in X_0$ be as in the statement of Proposition~(3.18).
By openness of $\mu$-stability the generic point of $X_0$ parametrizes a
$\mu$-stable sheaf.  Furthermore, since we have an  injection
$$\,\,H^0(F,F\ot K)^0\hra H^0(F^{**},F^{**}\ot K)^0\,,$$
we conclude by~(3.19) that $\cMx$ is good at $[F]$. This proves that $\cMx$
is good. The only thing to prove is that the generic point of $X$
parametrizes a locally-free sheaf. Since $F$ is stable, a neighborhood of
$[F]$ in $\cMx$ is isomorphic to $Def^0(F)$, which is smooth becuse $\cMx$
is good at $[F]$. Thus it suffices to show that there exists $x\in Def^0(F)$
parametrizing a locally-free sheaf. This is equivalent to the existence of
$x\in Def(F)$ parametrizing a locally-free sheaf. As is well-known this
follows from~(3.19); we will recall the proof. Let
$$\,\,{\rm Supp}\left(F^{**}/F\right)=\{P_1,\ldots,P_{\ell}\}\,.$$
The versal deformation space of $F_{P_i}$ (the localization of $F$ at
$P_i$) is smooth, and $F_{P_i}$ deforms to a free $\cO$-module~[F]. Hence,
since $Def(F)$ is smooth, it suffices to show that the map of tangent
spaces
$$T_0Def(F)\brel\rho\over\lra \oplus_{i=1}^{\ell}Def(F_{P_i})$$
induced by a versal sheaf on $S\tm Def(F)$, is surjective.
The map $\rho$ is part of the local-to-global exact sequence coming from
the spectral sequence abutting to ${\rm Ext}^1(F,F)$. The piece of interset
to us is
$$\,\,{\rm Ext}^1(F,F)\brel\rho\over\lra
H^0\left(\oplus_{i=1}^{\ell} Ext^1(F_{P_i},F_{P_i})\right)\to
H^2\left(Hom(F,F)\right)\,.\eqno(*)$$
We have an exact sequence
$$\,\,0\to Hom(F,F)\to Hom(F^{**},F^{**})\to R\to 0\,,$$
where $R$ is an Artinian sheaf (supprorted at the $P_i$'s). Hence~(3.19)
implies that the last term of~($*$) is zero, thus $\rho$ is surjective.
This completes the proof of Theorem~C.
\bsk

\n
{\bf Surfaces with ample canonical bundle.}
\msk
\n
Theorem~C  together with~(D2def) gives the following

\proclaim (3.26) Proposition.
Let $S$ be a surface with $K$ ample. Assume that there exists an effective
divisor $H$ on $S$ such that $c_1(K)=kc_1(H)$ for some rational positive
$k$, and such that $|n_0H|$ contains a smooth curve, where $n_0$ is given
by~(2.13). Let $N_2(S,H)$ be as in~(2.15). Let $\xi$ be a set of sheaf data
with $\rx=2$, and let $\cMx$ be the  corresponding moduli space of
sheaves on $S$,  polarized by $K$.  If
$$\,\,\Dx>N_2(S,H)K^2+2K^2+7\chi(\cO_S)\,,$$
then $\cMx$ is good.

\n
{\it Proof of Corollaries~$C'$ and~$C''$.}
\hskip 2mm
For Corollary~${\rm C}'$ notice that the hypotheses of Proposition~(2.11) are
satisfied by $H=K$. In fact $K$ is effective by hypothesis. Furthermore an
easy computation gives $z_0>2$, hence $n_0\ge 3$. Since $K^2$ is large
($K^2\ge 6$ suffices), $n_0K$ is very ample by~[Bo], in particular  there
exists a smooth curve in the linear system $|n_0K|$. The result follows
from Proposition~(3.26) and the easy estimate
$$\,\,N_2(S,K)<40+{22(\chi+1)\over 3K^2}\,,$$
valid for $K^2>100$. Similarly Corollary~${\rm C}''$
follows from Proposition~(3.26) and the estimate
$$\,\,N_2(S,H)<15+{8(\chi+1)\over 3K^2}\,,\eqno(3.27)$$
valid for $k>100$.
\msk

\n
{\it Examples of non-good moduli spaces.}
\hskip 2mm
We will show  that the lower bound given
in  Corollaries~${\rm C}'$-~${\rm C}''$ is, if not sharp, at least of the right
form.  Assume $H$ is
effective. We will consider non-trivial extensions
$$\,\,0\to\cO_S\to F\to I_Z(H)\to 0\,,\eqno(3.28)$$
where $Z$ is a zero-dimensional subscheme of $S$ of length $\ell$. These
extensions are parametrized by
$$\,\,{\rm Ext}^1\left(I_Z(H),\cO_S\right)\cong
H^1\left(I_Z(K+H)\right)^{*}\,.$$
{}From this one easily gets

\proclaim (3.29).
Keep notation as above. Assume that
$$\,\,\ell>h^0(K+H)=
\chi(\cO_S)+{1\over 2}K\cdot H+{1\over 2}H^2\,.\eqno(3.30)$$
Then if $Z$ is generic we have
$$\,\,\dim {\rm Ext}^1\left(I_Z(H),\cO_S\right)=
\ell-\chi(\cO_S)-{1\over 2}K\cdot H-{1\over 2}H^2\,,$$
and the generic non-trivial extension~(3.28) is locally-free.

\n
The following is also an easy exercise.

\proclaim (3.31).
Keeping notation as above, assume that~(3.30) is satisfied and that
$$\ell>{9\over 8}H^2+2+q_S\,.\eqno(ell2)$$
Then if $Z$ is generic the generic non-trivial extension~(3.28) is
$\mu$-stable, and furthermore its space of global sections is
one-dimensional.

\n
Let
$$\,\,\xi(\ell):=(2,\cO_S(H),\ell)\,.$$
By~(3.29)-(3.31) if
$$\ell>\max\left\{\chi(\cO_S)+{1\over 2}K\cdot H+{1\over 2}H^2,
{9\over 8}H^2+2+q_S\right\}\eqno(3.32)$$
then there is an irreducible subset $\S_{\ell}\ss\cM_{\xi(\ell)}$
parametrizing extensions~(3.28) with $Z$ generic, and we have
$$\,\,\dim\S_{\ell}=
3\ell-\chi(\cO_S)-{1\over 2}K\cdot H-{1\over 2}H^2\,.$$
Since $\D_{\xi(\ell)}=\ell-H^2/4$ we conclude that

\proclaim (3.33).
Let notation be as above. Assume that $\ell$ satisfies~(3.32) and that
$$\,\,\ell< 2\chi(\cO_S)-{1\over 2}K\cdot H+{1\over 2}H^2\,.$$
Then $\cM_{\xi(\ell)}$ is not good. More precisely it contains a subset
$\S_{\ell}$ whose dimension is greater than the expected dimension of
$\cM_{\xi(\ell)}$. Furthermore the generic point of $\S_{\ell}$ parametrizes
a $\mu$-stable vector bundle.

\n
One can easily check that if $q_S=0$ and if $c_1(K)=kc_1(H)$ for $k\gg 0$,
then the hypotheses of~(3.33) are satisfied. Hence letting
$$\,\,\ell_0:=2\chi(\cO_S)-{1\over 2}K\cdot H+{1\over 2}H^2-1\,,$$
one gets

\proclaim (3.34).
 Let notation be as above. Assume that $c_1(K)=kc_1(H)$ for $k\gg 0$. Then
\msk
\item{1.} The moduli space $\cM_{\xi(\ell_0)}$ is not good, and
\item{2.} $\D_{\xi(\ell_0)}=2\chi(\cO_S)-{1\over 2}K\cdot H+{1\over
4}H^2-1$.

\n
By considering surfaces with arbitrarily large $k$, we conclude that a
bound of the form given in  Corollaries~${\rm C}'$-~${\rm C}''$ is the best
we can hope for.
\bsk
\bsk

\n
{\xlbf 4. Irreducibility.}
\msk
\n
We will first derive Theorem~D  from
Theorem~A and Corollary~${\rm B}'$; the argument is due to
Gieseker-Li~[GL1]. Then we will  obtain Theorem~E  by
making explicit  Gieseker-Li's proof in the case considered by the
theorem.
\bsk

\n
{\bf Proof of Theorem~D.}
\msk
\n
We will prove  the following

\proclaim (4.1) Proposition.
Let $(S,H)$ be a polarized surface. For any integer $r\ge 2$ and any line
bundle $M$ on $S$, there exists a number $\D_3(r,M,S,H)$ such that the
following holds. Let $\xi$ be  set of sheaf data with
$\detx\cong M$, and such that
$$\,\,\Dx>\D_3(\rx,M,S,H)\,.\eqno(4.2)$$
Then $\cMx$ is irreducible (and is the closure of the subset parametrizing
$\mu$-stable vector bundle).

\n
This result implies Theorem~D. In fact, since tensorization by a line
bundle $N$ identifies $\cMx$ with $\cM_{\xi\ot N}$ (with the obvious
notation), Theorem~D will hold if we set
$$\,\,\D_3(r,S,H):=\max\left\{\D_3(r,M,S,H)\right\}_{M\in \cS}\,,$$
where $\cS$ is any  set of line bundles whose first Chern classes are a
set of representatives for the finite group
$H_{\ZZ}^{1,1}(S)/\rx H_{\ZZ}^{1,1}(S)$. The following lemma is the key
ingredient in the proof of Proposition~(4.1).

\proclaim (4.3) Lemma.
Let $(S,H)$ be a polarized surface. For any integer $r\ge 2$   there exists a
number $\wh{\D}_1(r,S,H)$ (with $\wh{\D}_1\ge\D_1'$, where $\D_1'$ is
as in Corollary~${\rm B}'$) such that the following holds. Let $\xi$, $\ell$
be a set of sheaf data and  a positive integer respectively, such that
$$\,\,\Dx\ge \wh{\D}_1(\rx,S,H)+\ell\,.$$
Set
$$\,\,\xi':=(\rx,\detx,c_2(\xi)-\ell)\,,$$
Let $X$ be an irreducible component of $\cMx$. Then there exist a locally
closed $\ul{\hb{non-empty}}$ subset $Y\ss \pal X$ and an open subset $V$
of an  irreducible component of $\cM_{\xi'}$,  with the following
properties. If $[E]\in V$ then $E$ is locally-free, $\mu$-stable, and
$$\,\,h^0(E,E\ot K)^0=0\,.$$
Furthermore $[F]\in Y$ if and only if $F$ fits
into an exact sequence
$$\,\,0\to F\to E\brel\phi\over \to
\oplus_{i=1}^{\ell}\CC_{P_i}\to 0\,,\eqno(4.4)$$
for some $[E]\in V$ (here $\{P_1,\ldots,P_{\ell}\}$ is a set of
distinct points of $S$).

\pf
First one proves the lemma  in the case $\ell=1$, then the general case
follows easily from this. The case $\ell=1$ follows from Theorem~A,
Corollary~${\rm B}'$ and dimension counts (use Proposition~(3.3) and
Theorem~(3.5)). The argument is similar to that used in the proof of
Corollary~${\rm B}'$; we leave the details to the reader.
\qed

\msk
Now we are ready to prove Proposition~(4.1). Set $r=\rx$.  Choose a set
of sheaf data
$$\,\,\xi_0=\left(r,M,c_2(\xi_0)\right)\,,$$
 such that
$$\,\D_{\xi_0}\ge\wh{\D}_1(\rx,S,H)
\hb{ with $c_2(\xi_0)$ minimal.}\eqno(4.5)$$
Notice that, since by~[HL,LQ] the moduli
space  $\cMx$ is non-empty for $\Dx\gg 0$,  Lemma~(4.3) shows that
also  $\cM_{\xi_0}$ is non-empty. Let
$U_{\xi_0}\ss\cM_{\xi_0}$ be the open subset parametrizing
$\mu$-stable locally-free sheaves $E$ such that
$$\,\,h^0(E,E\ot K)^0=0\,.$$
Since $\wh{\D}_1\ge\D_1'$, the subset $U_{\xi_0}$ is dense in
$\cM_{\xi_0}$, and hence non-empty. There exists an integer $n$ such
that for all $[E]\in U_{\xi_0}$ the bundle $E(nH)$ has $(r-1)$
independent sections, and such that the
degeneracy locus of the corresponding map
$$\cO_S^{(r-1)}\to E(nH)$$
is (at most) zero-dimensional. (For example if $E(nH)$ is
generated by  global sections.)
Hence if $[E]\in U_{\xi_0}$ then $E$   fits into an exact sequence
$$\,\,0\to \cO_S(-nH)^{(r-1)}\to E
\to I_Z\ot M\ot [(r-1)nH]\to 0\,,\eqno(4.6)$$
where $I_Z$ is the ideal sheaf of a zero-dimensional subscheme $Z$.  We
can, and  will, assume that
$$\,\,h^1\left(M\ot[rnH+K]\right)=0\,.\eqno(4.7)$$
Set
$$\,\,\D_3(r,M,S,H):=\wh{\D}_1(r,S,H)+h^0(M\ot[rnH+K])+1\,.\eqno(4.8)$$
We will prove that Proposition~(4.1) holds with this value of $\D_3$.
Thus we assume that $\Dx$ satisfies
Inequality~(4.2). Set
$$\,\,\ell:=\Dx-\wh{\D}_1\,.\eqno(\sharp)$$
Let $X$ be an irreducible component of $\cMx$. By
Lemma~(4.3) there exists $[E]\in U_{\xi_0}$ such that $X$ contains
all the isomorphism classes of sheaves $F$ fitting into~(4.4), with $E$
the chosen vector-bundle, and $\ell$ given by~($\sharp$). Since $E$ fits into
Exact Sequence~(4.6), if we choose $\phi$ appropiately we can arrange
that $F:=\ker\phi$ fit into the exact sequence
$$\,\,0\to \cO_S(-nH)^{(r-1)}\to F\to
 I_W\ot M\ot[(r-1)nH]\to 0\,,\eqno(*)$$
with $W=Z\cup\{P_1,\ldots,P_s\}$, for $s=\left(\Dx-\D_{\xi_0}\right)$.
Furthermore, since by hypotheses
$$\,\,\Dx-\D_{\xi_0}>h^0(M(rnH+K))\,,$$
we can also assume (using~(4.7)) that
$$\,\,h^1(I_W\ot M\ot[rnH+K])=
\ell(W)-h^0(M\ot[rnH+K])\,,$$
and hence by Serre duality
$$\,\,\dim{\rm Ext}^1
\left(I_W\ot M\ot [(r-1)nH],\cO_S(-nH)^{(r-1)}\right)=
(r-1)\cdot\left(\ell(W)-h^0(M\ot[rnH+K])\right)\,.\eqno(\dag)$$
Now let $\S_{\ell}$ be the  space  parametrizing non-trivial
extensions~($*$), such that~($\dag$) holds. Then $\S_{\ell}$ fibres
over a non-empty open subset of $Hilb^{\ell}(S)$, with projective spaces as
fibres. In particular $\S_{\ell}$ is irreducible.  Since $E$ is $\mu$-stable
we conclude that the subset  $\S_{\ell}^{\mu}\ss\S_{\ell}$ parametrizing
$\mu$-stable extensions is non-empty. Let $\O_{\ell}\ss\cMx$ be the
image of $\S_{\ell}^{\mu}$
under the classifying map; since $\S_{\ell}$ is irreducible, so is
$\O_{\ell}$. We have proved that $X$ contains a point
$[F]\in\O_{\ell}$.  Since $h^0(F^{**},F^{**}\ot K)^0$ vanishes, $\cMx$ is
smooth at $[F]$. Hence we conclude that $X$ contains all of
$\O_{\ell}$. To sum up: every irreducible
component of $\cMx$ contains the irreducible (non-empty) subset
$\O_{\ell}$, and  $\cMx$ is smooth at the generic point of
$\O_{\ell}$.  This implies that $\cMx$ is irreducible.

\proclaim (4.9) Remark.
The above proof works also if we only assume that for the generic
$[E]$ in any irreducible component of $U_{\xi_0}$,  $E$
fits into~(4.6).

\bsk

\n
{\bf Complete intersections with Picard number one.}
\msk
\n
The goal of this subsection is to prove the following

\proclaim (4.10) Proposition.
Let $(S,H)$ be as in the statement of Theorem~E. Assume also that
$Pic(S)=\ZZ[H]$. Then the conclusion of Theorem~E holds for
$(S,H)$.

We begin by giving an explicit value for $\wh{\D}_1$ for rank-two
sheaves.

\proclaim (4.11) Proposition.
Let $(S,H)$ be a polarized surface. Assume that $H$ is effective and that
the linear system $|n_0H|$ contains a smooth curve, where $n_0$ is
defined by~(2.13). Then Lemma~(4.3) holds (in rank two) with
$$\,\,\wh{\D}_1(2,S,H):=\D_2(S,H)+2h^0(2K)\,,$$
where $\D_2(S,H)$ is defined by~(2.12).

\pf
First we prove that~(4.3) holds for $\ell=1$. Thus we assume that
$\Dx\ge(\wh{\D}_1+1)$; in particular $S$, $H$ and $\xi$ satisfy the
hypotheses of Proposition~(3.18). Let $X$ be an irreducible component of
$\cMx$. Let $[F_0]\in \pal X^{\mu}$ be such that Proposition~(3.18) holds
with $F=F_0$. Let $\wt{Y}_1\ss\pal X^{\mu}$ be the irreducible component
of the stratum of the double dual stratification which contains $[F_0]$. Let
$Y_1\ss\wt{Y}_1$ be the open subset parametrizing sheaves $F$ such that
$$\,\,h^0(F^{**},F^{**}\ot K)^0=0\,.$$
The subset $Y_1$ is non-empty because $[F_0]\in Y_1$. We have
$Y_1^{**}\ss\cM_{\xi'}$, where
$$\,\,\xi'=(2,\detx,c_2(\xi'))\,,\quad c_2(\xi')<c_2(\xi)\,.\eqno(\sharp)$$
By Proposition~(3.3) and Theorem~(3.5)
$$\,\,\dim Y_1^{**}+3\left(c_2(\xi)-c_2(\xi')\right)\ge
\dim Y_1\ge\dim X-1\,.\eqno(*)$$
By~(3.19) both $\cMx$ and $\cM_{\xi'}$ are good at $[F_0]$ and
$[F_0^{**}]$ respectively. Hence~($*$) gives
$$\,\,\D_{\xi'}\ge\Dx-1\,.\eqno(\dag)$$
By~($\sharp$) we conclude that $c_2(\xi')=\left(c_2(\xi)-1\right)$.
Thus~($\dag$) is an equality, and hence also  all the inequalities
in~($*$). The result is that $\cM_{\xi'}$ is good (by Proposition~(3.18))
and furthermore
$$\,\,\dim Y_1^{**}=\dim\cM_{\xi'}\,.\eqno(\star)$$
Now set $V:=Y_1^{**}$. By~($\star$) $V$ is an open subset of $\cM_{\xi'}$.
By construction, if $[F]\in Y_1$ then $F$ fits into an exact
sequence~(4.4) for some $[E]\in V$ (with $\ell=1$). Conversely, since the
inequalities of~($*$) are equalities, if $[F]$ is the generic sheaf fitting
into~(4.4) for $[E]\in V$ (with $\ell=1$)  then $[F]\in Y_1$. We conclude that
$$Y:=\{[F]\in\cMx|
\hb{ $F$ fits into~(4.4) for some $[E]\in V$, with $\ell=1$}\}$$
is contained in the closure of $Y_1$, and hence $Y\ss X$. This proves
Lemma~(4.3) if $\ell=1$. When $\ell>1$ one iterates this construction.
Define $Y_1\ss\pal X$ as above. Let $X_2:=Y_1^{**}$, and define
$Y_2\ss\pal \ov{X}_2^{\mu}$ in the same way as we defined $Y_1\ss\pal
X^{\mu}$. We continue this process up to $Y_{\ell}$ (that this is possible is
guaranteed by Proposition~(3.18)).  Set $V:=Y_{\ell}^{**}$, and let
$Y\ss\cMx$ be the parameter space for all sheaves $F$ fitting
into~(4.4) for $[E]\in V$. One checks easily that $Y\ss\pal X$ and hence
that Lemma~(4.3) holds.
\qed

\proclaim (4.12) Corollary.
Let $(S,H)$ be a polarized surface satisfying the hypotheses of
Theorem~E. Then Lemma~(4.3) holds with
$$\,\,\wh{\D}_1(2,S,H)=17K^2+10\chi(\cO_S)\,.$$

\pf
Immediate from Proposition~(4.11) together with~(2.16) and~(3.27).
\qed

\msk
\n
Set
$$\,\,\xi_0:=(2,\cO_S,17K^2+10\chi(\cO_S))\,.$$
Then $\D_{\xi_0}$ satisfies~(4.5).

\proclaim (4.13) Lemma.
Keep notation as above. Let $(S,H)$ be a polarized surface satisfying the
hypotheses of Proposition~(4.10). Let $X\ss\cM_{\xi_0}$ be an irreducible
component. There is an open dense subset $U\ss X$, parametrizing
locally-free sheaves, such that if $[E]\in U$ then $E$ fits into an exact
sequence
$$\,\,0\to \cO_S(-6K)\to E\to I_Z\ot[6K]\to 0\,,\eqno(4.14)$$
where $Z$ is a zero-dimensional subscheme of $S$ of length
$$\,\,\ell(Z)=53K^2+10\chi(\cO_S)\,.$$

\pf
We will prove that $E$ fits into~(4.14); once this is done   the length of
$Z$ is obtained by computing  $c_2(E)$. We begin by showing
that if $[E]\in\cM_{\xi_0}$ and $E$ is locally-free, then $E\ot[5K]$ fits
into an exact sequence
$$\,\,0\to \cO_S(nH)\to E\ot[5K]\to I_W\ot[10K-nH]\to 0\,,\eqno(*)$$
where $W$ is a zero-dimensional subscheme of $S$, and  $n\ge 0$. In
fact by Serre duality and $\mu$-semistability we have
$$\,\,h^2(E\ot[5K])=h^0(E\ot[-4K])=0\,,$$
hence $h^0(E\ot[5K])\ge\chi(E\ot[5K])$. Applying the
Hirzebruch-Riemann-Roch Theorem one gets that
$$\,\,\chi(E\ot[5K])=3K^2-8\chi(\cO_S)\,.$$
The right-hand side is positive for $k\gg 0$, e.g.~if $k>100$, and hence
$h^0(E\ot[5K])>0$. (Here $k$ is as in the statement of Theorem~E.)
That $E$ fits into~($*$) follows from this and the hypothesis that
$Pic(S)=\ZZ[H]$ (in fact this is the only place where we use this
assumption). We need to bound  $n$. Of course by semistability we have
$n\le 5k$; we will show that if $[E]$ is generic there is a better bound.

\proclaim Claim.
Assume that $[E]\in X$ is generic. Then $n<4k$.

\n
{\bf Proof of the claim.}
\hskip 2mm
Assume that for generic $[E]\in X$  the sheaf $E\ot[5K]$ fits into~($*$)
 with $n\ge 4k$. Then
$$\,\,\dim\cM_{\xi_0}\left({K\cdot H\over \sqrt{H^2}}\right)
=4\D_{\xi_0}-3\chi(\cO_S)\,.$$
In fact this follows by writing
$$\,\,\mu(nH)=5K\cdot H-\a\sqrt{H^2}
\quad\hb{for}\quad\a\le {K\cdot H\over \sqrt{H^2}}\,.$$
Applying Proposition~(5.10) we get
$$\,\,3{(K\cdot H)^2\over H^2}+
{(K\cdot H)^2+2(K\cdot H)H^2+2K\cdot H\over H^2}+
{3\left((K\cdot H)+H^2+1\right)^2\over 2H^2}+
3\ge 17K^2+10\chi(\cO_S)\,.$$
As is easily verified this is false as soon as $k>0$, and hence we conclude
that $n<4k$ for generic $[E]\in X$ (with $E$ locally-free).
\qed

\msk
\n
 Now we are ready to show that if $[E]\in X$ is generic and $E$ is
locally-free, then $E\ot[6K]$ has a section with isolated zeroes. By the
above claim the vector-bundle  $E\ot[6K]$ fits into an exact sequence
$$\,\,0\to \cO_S(K+nH)\to E\ot[6K]\to
I_W\ot[11K-nH]\to 0\,,\eqno(\star)$$
with $n<4k$. First we show that
$$\,\,H^0\left(I_W\ot[11K-nH]\right)\not=0\,.$$
For this it suffices to prove that
$$\,\,h^0\left(\cO_S(11K-nH)\right)>\ell(W)\,.\eqno(\dag)$$
The left-hand side equals $\chi\left(\cO_S(11K-nH)\right)$, hence is
computed by Hirzebruch-Riemann-Roch.
An easy computation then gives that for~($\dag$) to hold we need that
$$\,\,(13k^2-2k-2)H^2>9\chi(\cO_S)\,.$$
 This inequality is satisfied as soon as $k>2$. Let
$$\s\in H^0\left(I_W\ot[11K-nH]\right)$$
be a non-zero section; it lifts to section
$\wt{\s}$ of $E\ot[6K]$ because $h^1(K+nH)=0$. If $\tau$ is any section of
$\cO_S(K+nH)$, then
$$\,\,(\tau+\wt{\s})\ss (\s)\,,$$
where $(\cdot)$ denotes "zero-locus". Since $\cO_S(K+nH)$ is very
ample one easily concludes that there exists $\tau$ such that
$\t:=(\tau+\wt{\s})$ is section with isolated zeroes. Thus $E$ fits into
$$\,\,0\to \cO_S\brel\t\over\lra E\ot[6K]\to I_Z\ot[12K]\to 0\,,$$
where $Z=(\t)$. Tensoring the above exact sequence
with $\cO_S(-6K)$ one obtains~(4.14).
\qed
\msk

\n
{\it Proof of Proposition~(4.10).}
\hskip 2mm
By Remark~(4.9), Formula~(4.8), Lemma~(4.13) and Corollary~(4.12),
we can set
$$\,\,\D_3(2,\cO_S,S,H)=\wh{\D}_1+h^0(13K)+1\,,$$
where the value of $\wh{\D}_1$ is given  by Corollary~(4.12).
Proposition~(4.10) follows at once.
\bsk

\n
{\bf Proof of Theorem~E.}
\msk
\n
Let $S$ be a surface satisfying the hypotheses of Theorem~E.
Then
$$\,\,S=V_1\cap\cdots\cap V_n\ss\PP^{n+2}\,,$$
where $V_i$ is a degree-$d_i$ hypersurface. Let
$$\rho\cl \cS\to B$$
 be the family of smooth complete intersections of $n$ hypersurfaces in
$\PP^{n+2}$ of degrees $d_1\ldots,d_n$. Thus $B$ is an open subset of a
Grassmannian, and $S=\rho^{-1}(b_0)$ for a certain $b_0\in B$. If $b\in B$
we set $S_b:=\rho^{-1}(b)$. We let $B_1\ss B$ be the subset parametrizing
surfaces whose Picard group is generated the hyperplane class. We recall
the following

\proclaim (4.15) Noether-Lefscetz Theorem.
Keep notation as above. Assume that $p_g(S_b)>0$  for $b\in B$. Then
$B_1$ is dense in $B$.

For $c\in\ZZ$ let $\cM_c(S_b)$ be the moduli space of
torsion-free sheaves $F$ on $S_b$, semistable with respect to
$\cO_{S_b}(1)$, with
$$\,\,r(F)=2\,,\quad \det F\cong\cO_{S_b}\,,\quad c_2(F)=c\,.$$
By Maruyama~[Ma] there exists a relative moduli space
$$\pi\cl\cM_c(\cS)\to B$$
proper over $B$, such that $\pi^{-1}(b)\cong\cM_c(S_b)$ for all $b\in B$.

\proclaim (4.16) Proposition.
Keep notation as above. Assume that the integer $k$ such that
$K_{S_b}\sim kH$ is large (e.g.~$k>100$). If
$$c>95K_{S_b}^2+11\chi(\cO_{S_b})+1$$
then $\cM_c(\cS)$ is irreducible.

\pf
For $b\in B$ let $\cM_c^0(S_b)\ss\cM_c(S_b)$ be the (open) subset
parametrizing sheaves $F$ such that
\msk
\item{1.} $F$ is locally-free and stable,
\item{2.} $h^0(F,F\ot K_{S_b})^0=0$.
\msk
\n
Let $\cM_c^0(\cS):=\cup_{b\in B}\cM_c^0(S_b)$. By Corollary~${\rm C}''$
$\cM_c^0(S_b)$ is dense in $\cM_c(S_b)$ for all $b$, and  hence
$\cM_c^0(\cS)$ is dense in $\cM_c(\cS)$. Thus it
suffices to show that $\cM_c^0(\cS)$ is irreducible. Let $X$ be anyone  of
its irreducible components. Let $X^{*}\ss X$ be the complement of the
intersection with all other irreducible components. We claim that
$\pi(X^{*})$ contains an open non-empty subset of $B$. Since $\pi(X^{*})$
is constructible it suffices to prove that it is not contained in any
proper subvariety of $B$.   Let $b\in\pi(X^{*})$, and let $v\in T_b(B)$. By
deformation theory (see for example Proposition 2.1 in~[G2]) there
exist a curve $\i\cl \L\hra X^{*}$ and a point $P\in\L$ such that
$\i(P)=[F]$ and
$$\,\,v\in\Im D(\pi\circ\i)(P)\,,$$
where $D$ is the differential. This proves that $\pi(X^{*})$ is not
contained in any proper subvariety of $B$, and hence it contains an
open non-empty subset. Now let $Y$ be another irreducible component of
$\cM_c^0(\cS)$. Then
$$\,\,X^{*}\cap Y^{*}=\es\,.\eqno(\star)$$
Since $\pi(X^{*})$ and $\pi(Y^{*})$ both contain a Zariski-open non-empty
subset of the irreducible variety $B$, we conclude by~(4.15) that there
exists
$$b_0\in\pi(X^{*})\cap\pi(Y^{*})$$
  such that $Pic(S_{b_0})\cong\ZZ[\cO_{S_{b_0}}]$. By
Proposition~(4.10) the moduli space $\cM_c(S_{b_0})$ is irreducible.
Since both $X^{*}$ and $Y^{*}$ must contain an open non-empty subset of
$\cM_c(S_{b_0})$, we conclude that
$$\,\,X^{*}\cap Y^{*}\not=\es\,.$$
This contradicts~($\star$), and  hence $\cM_c^0(\cS)$ is irreducible.
\qed

\proclaim (4.17) Corollary.
Keep assumptions as in Proposition~(4.16). Then $\cM_c(S_b)$ is
connected for all $b\in B$.

\pf
Let $b\in B_1$. By Proposition~(4.10), $\pi^{-1}(b)=\cM_c(S_b)$ is
irreducible, hence connected. Since $B_1$ is dense in $B$ (by~(4.15)), and
since $\cM_c(\cS)$ is irreducible (Proposition~(4.16)) and proper over $B$,
we conclude that $\cM_c(S_b)$ is connected for all $b\in B$.
\qed

\msk
Now we are ready to prove Theorem~E. Let $S$ and $H$ be as in
the statement of the theorem, and let $c$ be an integer such that
$$\,\,c>95K^2+\chi(\cO_S)+1\,.\eqno(cgr)$$
We must prove that $\cM_c(S)$ is irreducible.

\proclaim Claim.
Keep notation and assumptions as above.  Suppose also that $\cM_c(S)$
is reducible. Then there exist two irreducible components $X_1$, $X_2$
such that $X_1\cap X_2$ contains a point  parametrizing a  stable
sheaf.

\pf
By Corollary~(4.17) there exist two irreducible components $X_1$,
$X_2$ such that their intersection is non-empty. Let $[F]\in X_1\cap X_2$.
If $F$ is stable there is nothing to prove, so assume $F$ is non-stable. Let
$E:=Gr(F)$. By Luna's \'etale slice Theorem the natural map
$$\l\cl Def(E)\to \cM_c(S)$$
surjects onto a neighborhood of $[F]$. (Since $S$ is regular,
$Def^0(E)=Def(E)$.)  Hence $\l^{-1}X_i$ is a closed non-empty subset of
$Def(E)$, and
$$\,\,\dim\l^{-1}X_i\ge\dim X_i=4c-3\chi(\cO_S)\,.\eqno(*)$$
On the other hand the dimension of the tangent space to $Def(E)$ at the
origin is given by
$$\,\,\dim T_0\left(Def^0(E)\right)=
h^1(E,E)=-\chi(E,E)+h^0(E,E)+h^0(E,E\ot K)\,.$$
Here $\chi(E,E)$ and the other terms are  defined as in~(5.1). By
Lemma~(5.41) and by~(5.7) we conclude that
$$\,\,\dim T_0Def(E)\le 4c-3\chi(\cO_S)+3+
{2\over H^2}(K\cdot H+5H^2+1)^2\,.$$
Thus, by Inequality~($*$), we have
$$\dim\left(\l^{-1}X_1\cap\l^{-1}X_2\right)\ge 4c-3\chi(\cO_S)
-3-{2\over H^2}(K\cdot H+5H^2+1)^2\,.$$
As is easily checked
$$\,\,4c-3\chi(\cO_S)-3-{2\over H^2}(K\cdot H+5H^2+1)^2>
3c+\e(2,S,H)\,,$$
if $k>0$. Hence by Proposition~(5.40) there exists a point in
$x\in\l^{-1}X_1\cap\l^{-1}X_2$ parametrizing a $\mu$-stable sheaf. Then
$\l(x)\in X_1\cap X_2$ parametrizes a stable sheaf.
\qed

\msk
\n
So let's assume that $\cM_c(S)$ is reducible. By the above claim  there
exist two irreducible components $X$, $Y$, of $\cM_c(S)$ and a point $[F]$
in their intersection such that $F$ is stable. By~(0.3) we conclude that
$$\,\,\dim_{[F]}X\cap Y\ge 4c-3\chi(\cO_S)-h^0(F,F\ot K)^0\,.$$
Applying~(5.7) we get
$$\,\,\dim_{[F]}X\cap Y\ge
4c-3\chi(\cO_S)-{2\over H^2}(K\cdot H+5H^2+1)^2\,.\eqno(\dag)$$
Since $X\cap Y$ is in the singular locus of $\cM_c(S)$ and since $F$ is
stable, we have
$$\,\,\dim_{[F]}X\cap Y\le W_{\xi}^K\,,$$
where $\xi=(2,\cO_S,c)$. By Theorem~B and by~($\dag$) we get
$$\,\,4c-3\chi(\cO_S)-{2\over H^2}(K\cdot H+5H^2+1)^2\le
\l_2c+\l_1\sqrt{c}+\l_0'+e_K\,.$$
A straightforward computation shows that this is impossible if
$c$ satisfies~(cgr) and $k$ is large, for example $k>100$. This proves
that $\cM_c(S)$ is irreducible.
\bsk
\bsk

\n
{\xlbf 5. Estimates.}
\msk
\n
In this section we prove various technical results which have been used in the
course of the paper. We let $H$ be an ample divisor on  the surface $S$.
Unless otherwise stated stablity and $\mu$-stability of sheaves on $S$ is
with respect to $H$. If $F$, $G$ are sheaves on the
same scheme we set
$$\,\,h^i(F,G):=\dim Ext^i(F,G)\,,\quad
\chi(F,G):=\sum_i(-1)^ih^i(F,G)\,.\quad(5.1)$$
\bsk

\n
{\bf A bound for the number of sections of a semistable sheaf on $S$.}
\msk
\n
The  bound provided by the following proposition is due to Simpson and Le
Potier~[S,Le].

\proclaim (5.2) Proposition (Simpson-Le Potier).
Keeping notation  as above, assume that $H$ satisfies~(0.4). Let $G$ be a
$\mu$-semistable torsion-free sheaf on $S$. Then
$$\,\,h^0(G)\le{r_G\over 2H^2}(\mu_G+(r_G+1)H^2+1)^2\,.$$

\n
The following corollary  is what we use.

\proclaim (5.3) Corollary.
Assume that $H$ satisfies~(0.4). Let $F,G$ be $\mu$-semistable
torsion-free sheaves on $S$. Then
$$\,\,h^0(F,G)\le{r_F\cdot r_G\over 2H^2}
(\mu_G-\mu_F+(r_F\cdot r_G+1)H^2+1)^2\,.$$

\pf
Since $F$ and $G$ are both torsion-free $\mu$-semistable sheaves on $S$,  so is
$Hom(F,G)$. The corollary follows by applying Proposition~(5.2) to $Hom(F,G)$.
\qed

\msk
In  the rank-one case there is a slightly sharper version of
Proposition~(5.2). If $A$ is an effective divisor on $S$, let
$$\,\,A_f:=\hb{fixed part of $|A|$,}\quad
A_m:=A-A_f=\hb{moving part of $A$}\,.$$

\proclaim (5.4) Lemma.
Let notation be as above. Assume that $H$ is effective. Let $C\in |H|$. If
$A$ is a divisor on $S$, then
$$\,\,h^0\left(\cO_C(A_m)\right)\le \left[C\cdot A +1\right]_{+}\,,$$
where $[x]_{+}:=\max\{x,0\}$.

\pf
We might as well assume that $A_m$ is effective
non-zero. Since $|A_m|$ has no fixed part, there exists $D\in |A_m|$
such that the  scheme-theoretic intersection of $C$ and $D$, call it $Z$, is
zero-dimensional.  Then the exact sequence
$$0\to \cO_C\to \cO_C(D)\to \cO_Z\to 0$$
gives
$$\,\,h^0\left(\cO_C(A_m)\right)\le h^0(\cO_C)+C\cdot D\le
h^0(\cO_C)+C\cdot A\,.$$
As is easily checked $C$ is 1-connected, hence $h^0(\cO_C)=1$ (see~[BPV]).
Hence the result follows from the above inequality.
\qed

\proclaim (5.5) Proposition.
 Assume that $H$ is effective.
Let $A$ be a rank-one torsion-free sheaf on $S$. Then
$$\,\,h^0(A)\le{1\over 2H^2}\left(\mu_A+H^2+1\right)^2\,.$$

\pf
Let $C\in|H|$. By considering the exact sequence
$$0\to \cO_S\left(A-(i+1)H\right)\to \cO_S(A-iH)\to \cO_C(A-iH)\to 0$$
for all non-negative integers $i$, one easily concludes that
$$\,\,h^0\left(\cO_S(A)\right)\le
\sum_{i=0}^{N}h^0\left(\cO_C((A-iH)_m)\right)\,,$$
where $N$ is the maximum value of $i$ such that
$$\,\,h^0\left(\cO_S(A-iH)\right)>0\,.$$
Thus by Lemma~(5.4)
$$\,\,h^0\left(\cO_S(A)\right)\le
\sum_{i=0}^{N}[\mu_A-iH^2+1]_{+}\,.$$
An easy estimate of the right-hand side gives the proposition.
\qed
\bsk

\n
{\bf Twisted endomorphisms.}
\msk
\n
If $L$ is a line bundle on $S$, and $r$ is a positive integer we set
$$\,\,e_L(r,S,H):=\cases{
{r^2\over 2H^2}\left(\mu_L+(r^2+1)H^2+1\right)^2-h^0(L)\,, &
if $L\cdot H\ge 0$,\cr
0\,, & if $L\cdot H<0$.}\eqno(5.6)$$
By applying Corollary~(5.3)  with $G=F\ot L$
one gets

\proclaim (5.7).
Assume that $H$ satisfies~(0.4). Let $F$ be a $\mu$-semistable rank-$r$
torsion-free sheaf on $S$. Then
$$\,\,\dim Hom(F,F\ot L)^0\le  e_L(r,S,H)\,.$$

\n
The following  is a useful observation:

\proclaim  Proposition.
Let $L$ be a line bundle on $S$. Let $F$ be a rank-two torsion-free sheaf on
$S$ such that
$$\,\,h^0(F,F\ot L)^0>h^0(L^{\ot 2})\,.$$
Then $F$ is $\ul{\rm not}$ $\left(L\cdot H/2\sqrt{H^2}\right)$-stable,
i.e.~there exists a rank-one subsheaf $A\ss F$ such that
$$\,\,\mu_A\ge \mu_F-{L\cdot H\over 2}\,.$$

\pf
Replacing $F$ by $F^{**}$ we can assume that  $F$ is locally-free.
Let
$$\phi\cl H^0(F,F\ot L)^0\to H^0(L^{\ot 2})$$
be the map defined by  $\phi(\s):=\det\s$. Since
$0\in\phi^{-1}(0)$, $\phi^{-1}(0)$ is non-empty,  and hence  by
our hypothesis $\dim\phi^{-1}(0)>0$.  Thus there exists a
non-zero map
$$\ \,f\cl F\to F\ot L\ \,$$
with $\det f={\rm tr}f=0$, i.e.~$f\circ f=0$.  The kernel of
$f$ is a rank-one torsion-free subsheaf of $F$, which we denote by $A$. We have
$$\,\,0\to A\to F\brel q\over \to B\to 0\,.$$
The sheaf $B$ is also torsion-free, because it is  isomorphic to $\Im f\ss F\ot
L$.
Since $f\circ f=0$, the map $f$  is obtained as the composition
$$\,\,F \brel q\over\to B\brel g\over \to A\ot L\,,$$
for some non-zero map $g$. Thus
$$\,\,\mu_B-\mu_A\le \mu_L\,.$$
This  implies the proposition.
\qed

\proclaim (5.8) Corollary.
Let $\xi$ be a set of sheaf data for $S$, with $\rx=2$. Let $\S_{\xi}\ss\cMx$
be the subset parametrizing  sheaves $F$ such that
$$\,\,h^0(F,F\ot K)^0>h^0(2K)\,.$$
Then $\S_{\xi}\ss\cMx\left((K\cdot H)/2\sqrt{H^2}\right)$.

\bsk

\n
{\bf A bound for the dimension of $\cM_{\xi}(\a)$.}
\msk
\n
Let
$$\eqalign{s(r):= & \left((r-1)^2+1\right)H^2+1\,,\cr
M(r,S,H):= &\cases{-{r\over 8}K^2-{1\over 2}r(r+1)\chi(\cO_S) &
if $\chi(\cO_S)\ge 0$ and $K^2\ge 0$,\cr
{r^2\over 2}|\chi(\cO_S)|+{r^3\over 8}|K^2| & otherwise.}}$$
The goal of this subsection is to prove the following  propositions.

\proclaim (5.9) Proposition.
Assume that $H$ satisfies~(0.4). Let $\xi$ be a set of sheaf data, and
$\a\in\RR$. Then $\cMx(\a)$ is a constructible subset of $\cMx$. If $K\cdot
H\ge 0$,
$$\eqalign{\dim\cM_{\xi}(\a)\le & (2r_{\xi}-1)\Dx+(2\rx-1)\a^2+
(\rx^2-2\rx+2){K\cdot H+2s(\rx)\over2\sqrt{H^2}}\a\cr
& +(\rx-1)(\rx^3-4\rx^2+6\rx-2){(K\cdot H)^2\over 8H^2}
+\left[(\rx-1)^3+\rx\right]{s(\rx)^2\over 2H^2}\cr
&+(\rx^2-\rx+1){(K\cdot H+s(\rx))^2\over 2H^2}
-q_S+M(\rx,S,H)\,.}$$
If $K\cdot H<0$,
$$\eqalign{\dim\cM_{\xi}(\a)\le & (2\rx-1)\Dx+(2\rx-1)\a^2+
\left[(\rx^2-2\rx+2){s(\rx)\over\sqrt{H^2}}-(\rx^2-2\rx-2)
{K\cdot H\over2\sqrt{H^2}}\right]\a\cr
& +{\rx^3\over 16}(17\rx+4){(K\cdot H)^2\over 8H^2}
+\left[(\rx-1)^3+\rx\right]{s(\rx)^2\over 2H^2}\cr
&+(\rx^2-\rx+1){(K\cdot H+s(\rx))^2\over 2H^2}
+{\rx^2\over 2}\left(1-2\chi(\cO_S)\right)+{\rx^3\over 8}|K^2|-q_S\,.}$$

\n
The next proposition provides a better bound for the case of rank two.

\proclaim (5.10) Proposition.
Assume that $H$ is effective.  Let $\xi$ be a set of sheaf data
with $r_{\xi}=2$, and let $\a\in\RR$. Then
$$\eqalign{\dim\cM_{\xi}(\a)\le &3\D_{\xi}+3\a^2
+{K\cdot H+2H^2+2\over\sqrt{H^2}}\a \cr
& +{3(K\cdot H+H^2+1)^2\over 2H^2}+{(K\cdot H)^2\over 4H^2}
-{K^2\over 4}+3-3\chi(\cO_S)-q_S\,.}$$

\n
The proposition below gives  a bound  for the dimension of the subset
$\cMx^C(\a)\ss\cMx(\a)$ defined in~(0.7).

\proclaim (5.11) Proposition.
Assume $H$ is effective. Let $C\in |nH|$ be a  smooth curve, and let $\xi$ be a
set of sheaf data with $\rx=2$. Let $\vf(\xi,\a)$ be the right-hand side of
the inequality in the previous proposition. Then $\cMx^C(\a)$ is a
constructible subset of  $\cMx$, and
$$\,\,\dim\cMx^C(\a)\le
\max\{\vf(\xi,\a_0)-n\a_0\sqrt{H^2}\}_{0\le\a_0\le\a}\,.$$ %

First we will show that the above propositions follow from a bound
for the dimension of certain extension groups (Propositions~(5.15) and~(5.16)).
Then we
will obtain these bounds.
\msk

\n
{\it Filtrations.}
\hskip 2mm
We will prove Propositions~(5.9)-(5.10) by bounding the number
of moduli of certain filtrations.  Let $F$
be a torsion-free sheaf on $S$. Let
$$0=F_0\ss F_1\ss F_2\ss\cdots\ss F_{n+1}=F\eqno(5.12)$$
be a filtration with $n\ge 1$. We set
$$\,\,Q_i:=F_i/F_{i-1}\,,\  r_i:= r(Q_i)\,,\
\mu_i:=\mu(Q_i)\,,\ \D_i:=\D(Q_i)\,.\eqno(5.13)$$
We also set $A:= F_1$, $B:= F/A$.
We will  make the following assumption

\proclaim (5.14).
\msk
\item{1.} $Q_i$ is torsion-free and $\mu$-semistable, and
\item{2.} $\mu_2>\mu_3>\cdots>\mu_{n+1}$.

\n
The significance of the following two propositions is that they  give an
upper bound on the number of moduli of such filtrations.

\proclaim (5.15) Proposition.
Keep notation  as above. Assume that $H$ satisfies~(0.4). Suppose that $F$
is $\mu$-semistable and that~(5.14) holds. Define $\a$ by setting
$$\,\,\mu_A=\mu_F-{\a\over r_A}\sqrt{H^2}\,.$$
If $K\cdot H\ge 0$, then
$$\eqalign{\sum_{i\le j}h^1(Q_j,Q_i)\le &  (2r_F-1)\D_F+
{1\over 2}(1+{1\over r_A})(2r_F-r_A)\a^2
+\left[r_B(r_B-1)+r_F\right]{K\cdot H+2s(r_F)\over2\sqrt{H^2}}\a\cr
&+r_B(r_B^3-r_B^2-r_Fr_B+r_F^2){(K\cdot H)^2\over 8H^2}
+(r_B^3+r_Ar_F){s(r_F)^2\over 2H^2}\cr
&+(r_F^2-r_Ar_B){(K\cdot H+s(r_F))^2\over 2H^2}
+M(r_F,S,H)\,.}$$
If $K\cdot H<0$, then
$$\eqalign{\sum_{i\le j}h^1(Q_j,Q_i) \le &  (2r_F-1)\D_F+
{1\over 2}(1+{1\over r_A})(2r_F-r_A)\a^2\cr
&+\left\{\left[r_F+r_B(r_B-1)\right]{K\cdot
H+2s(r_F)\over2\sqrt{H^2}}   -r_F(r_B-1){K\cdot
H\over\sqrt{H^2}}\right\}\a\cr
&+\left[r_Fr_B(r_B-1)(r_F+2r_A)+r_Ar_B^2(r_A+1)\right] {(K\cdot
H)^2\over 8H^2}+(r_B^3+r_Ar_F){s(r_F)^2\over 2H^2}\cr
&+(r_F^2-r_Ar_B){(K\cdot H+s(r_F))^2\over 2H^2}
+{r_F^2\over 2}\left(1-2\chi(\cO_S)\right)
+{r_F^3\over 8}|K^2|\,.}$$

\proclaim (5.16) Proposition.
Keeping notation as above, assume that $H$ is
effective. Suppose that  $F$ is
a rank-two $\mu$-semistable sheaf with a filtration~(5.12), and that~(5.14) is
satisfied.  Define $\a$ by setting
$$\,\,\mu_A=\mu_F-{\a\over r_A}\sqrt{H^2}\,.$$
Then
$$\eqalign{\sum_{i\le j}h^1(Q_j,Q_i)  \le & 3\D_F+3\a^2
+{K\cdot H+2H^2+2\over\sqrt{H^2}}\a \cr
&+{3(K\cdot H+H^2+1)^2\over 2H^2}+
{(K\cdot H)^2\over 4H^2}-{K^2\over 4}+3-3\chi(\cO_S)\,.}$$

\msk

\n
{\it Proof  of~(5.9)-(5.10)-(5.11) assuming~(5.15)-(5.16).}
\hskip 2mm
First we prove Propositions~(5.9) and~(5.10).  Let $\cF$ be a  family
of sheaves on $S$ parametrized by a scheme $B$. Let
$\ov{\chi}:=\left(\chi_1,\ldots,\chi_{n+1}\right)$ be a sequence of
one-variable polynomials. Drezet-Le Potier~[DL] have constructed a scheme
$$\,\,Drap(\cF;\ov{\chi})\to B\,,$$
proper over $B$, parametrizing all filtrations~(5.12) with
$F=\cF_x$ for some $x\in B$, and where
$$\,\,\chi\left(Q_i(nH)\right)=\chi_i\quad
i=1,\ldots,(n+1)\,.\eqno(\sharp)$$
Now let $\Px$ be a parameter
space for semistable sheaves of which $\cMx$ is the geometric invariant
theory quotient. Thus there is a family of semistable sheaves $\cFx$ on $S$
parametrized by $\Px$ inducing a surjection $\Px\to\cMx$.     For
$\a_0\in\RR$ let $I(\a_0)$ be the set of sequences $\ov{\chi}$ such that if
$F$ fits into~(5.12) and ~($\sharp$) holds, then
$$\,\,\mu(A)\ge \mu_F-{\a_0\over r(A)}\sqrt{H^2}\,.\eqno(5.17)$$
(The definition of $I(\a_0)$ makes sense because the slope of a sheaf is
determined by its Hilbert polynomial.) If
$\ov{\chi}$ is as above, let $Drap_0(\cFx;\ov{\chi})$ be the subset of
$Drap(\cFx;\ov{\chi})$ parametrizing filtrations such that~(5.14) is
satisfied. Since each of the properties of~(5.14) is open,
$Drap_0(\cFx;\ov{\chi})$ is an open subset of $Drap(\cFx;\ov{\chi})$. Set
$$\,\,Drap_0(\cFx,\a_0):=
\bigcup_{\ov{\chi}\in I(\a_0)}Drap_0(\cFx;\ov{\chi})\,.$$

\proclaim Claim.
Let notation be as above. Then $Drap_0(\cFx,\a_0)$ is a scheme of finite
type.

\pf
Since $Drap_0(\cFx;\ov{\chi})$ is of finite type for every $\ov{\chi}$, all we
have to show is that
$$I^{*}(\a_0):=
\left\{\ov{\chi}\in I(\a_0)|\ Drap_0(\cFx;\ov{\chi})\not=\es\right\}$$
is a finite set.
So let  $F$ be a sheaf satisfying the hypotheses of Proposition~(5.15),
and let $Q_i$ be as in~(5.13).  Since the set of  slopes of subsheaves
$Q_1\ss\left(\cFx\right)_x$, for $x\in\Px$, satisfying~(5.17) is  finite, it
suffices  to check that there is a universal bound for the size of the
coefficients of the Hilbert polynomials of the $Q_i$, if $(S,H)$, $\a$, $r_F$,
$c_1(F)$ and $c_2(F)$ are fixed. This in turn amounts to bounding the size of
$$\,\,c_1(Q_i)\cdot H\,,\quad c_1(Q_i)^2\,,
\quad c_1(Q_i)\cdot K\,,\quad c_2(Q_i)\,.$$
That  $c_1(Q_i)\cdot H$ is bounded follows from~(5.17) and from Item~(2)
of~(5.14). Thus, by Hodge index, we also get that $c_1(Q_i)^2$ is bounded
above. Let's show that it is also bounded below. By Item~(1) of~(5.14) and
Bogomolov's Inequality we have
$$\,\,{1\over 2r_i}c_1(Q_i)^2\ge {1\over 2}c_1(Q_i)^2-c_2(Q_i)\,.$$
Thus
$$\,\,\sum_{i=1}^{n+1}{1\over 2r_i}c_1(Q_i)^2\ge
\sum_{i=1}^{n+1}{1\over 2}c_1(Q_i)^2-c_2(Q_i)=
{1\over 2}c_1(F)^2-c_2(F)\,.\eqno(\star)$$
Since the values of $c_1(Q_i)^2$ are bounded above one concludes that they
are also bounded below. Since $c_1(Q_i)\cdot H$ is bounded we conclude
by the Hodge index theorem  that $c_1(Q_i)\cdot K$ is
bounded. Finally boundedness of $c_2(Q_i)$ follows from boundedness of
$c_1(Q_i)^2$, Bogomolov's Inequality and the equality in~($\star$).
\qed

\msk
\n
Let $\pi$ be the composition of the projection
$Drap_0(\cFx,\a_0)\to\Px$ and the quotient map $\Px\to\cMx$. By
construction we have
$$\,\,\cMx(\a_0)=\pi\left(Drap_0(\cFx,\a_0)\right)\,.\eqno(\flat)$$
Hence, by the claim above, $\cMx(\a_0)$
is a constructible set. Now let's consider $\dim\cMx(\a_0)$. For
convenience of exposition we will assume that there is a tautological
family of sheaves $\cGx$ parametrized by $\cMx$; it will be clear how to
modify the argument if this is not the case.   Let $x\in
Drap\left(\cGx;\ov{\chi}\right)$ correspond to~(5.12).
An easy inductive argument with extensions shows that
$$\,\,\dim Drap\left(\cGx;\ov{\chi}\right)_x\le
\sum_{i\le j}h^1(Q_j,Q_i)\,.$$
By~($\flat$) we conclude that $\dim\cMx(\a_0)$ is bounded above by the
maximum of the values of the right-hand side of the inequality in~(5.15)
(or in~(5.16)) for $r_F=\rx$, $\D_F=\Dx$,
$0\le\a\le\a_0$, and $1\le r_A\le \rx$. This immediately gives~(5.10). It
also gives a slightly weaker version of~(5.9). In order to get~(5.9) one
argues that if $[F]\in\cMx(\a_0)$ then at least one of $F$, $F^{**}$ fits into
a
filtration~(5.12) such that~(5.14) is satisfied and furthermore $r(A)\le
\rx/2$. The arguments above together with Theorem~(3.5) will show then
that  $\dim\cMx(\a_0)$ is bounded above by the
maximum of the values of the right-hand side of the inequality in~(5.15)
for $r_F=\rx$, $\D_F=\Dx$,
$0\le\a\le\a_0$, and $1\le r_A\le \rx/2$. Proposition~(5.9) follows
 from this together with easy estimates.

Now let's  prove  Proposition~(5.11).
We keep  the notation introduced in the previous proof. For
simplicity of exposition we assume that there exists a tautological family
$\cGx$ of sheaves parametrized by $\cMx$. Let $Drap_0^C(\cGx;\a_0)$ be the
subset of $Drap_0(\cGx;\a_0)$ parametrizing filtrations
$$\,\,0\to A\brel f\over\lra F\to B\to 0\,,\eqno(5.18)$$
such that the restriction $A|_C$
is locally-free and  spans a destabilizing subline bundle of $F|_C$.
As is easily checked $Drap_0^C(\cGx;\a_0)$ is closed. Since
$$\,\,\pi\left(Drap_0^C(\cGx;\a_0)\right)=\cMx^C(\a_0)\,,$$
we conclude that $\cMx^C(\a_0)$ is constructible. Now let's prove the upper
 bound for its dimension. We start by examining the map of vector bundles
$f|_C$. Let $\O$ be its zero-locus. Since $A|_C$ spans a
destabilizing subline bundle of $F|_C$, we have
$$\,\,\mu\left(A|_C\right)+\deg\O\ge \mu\left(F|_C\right)\,.\eqno(*)$$
Define $\a$ by setting
$$\,\,\mu_A=\mu_F-\a\sqrt{H^2}\,.$$
(Thus $0\le\a\le\a_0$.) Since $C\in |nH|$ we conclude by~($*$)  that
$$\,\,\deg\O\ge n\a\sqrt{H^2}\,.\eqno(5.19)$$
Since $B$ is a rank-one torsion-free sheaf we have  $B=I_Z\ot B^{**}$, for a
zero-dimensional subscheme $Z$ of $S$. We claim that for each point $P\in \O$
we have
$$\,\,{\rm mult}_P(Z)\ge {\rm mult}_P(\O) \,.\eqno(5.20)$$
In fact, let $f=(f_1,f_2)$ be an expression for $f$ in a neighborhood of $P$,
so that  $I_Z$ is
locally generated by $f_1,f_2$. Let $y$ be a local equation for $C$.  Then
$$\,\,{\rm mult}_P(Z)=\dim_{\CC}\cO/(f_1,f_2),
\quad {\rm mult}_P(\O)=\dim_{\CC}\cO/(f_1,f_2,y)\,,$$
where $\cO$ is the local ring at $P$; Inequality~(5.20) follows at once
from this.  Putting together~(5.20) and~(5.19) we get
$$\,\,\sum_{P\in C}{\rm mult}_P(Z)\ge n\a\sqrt{H^2}\,.\eqno(5.21)$$
To conclude the proof of Proposition~(5.11) we  go back to the proof
that~(5.16) implies~(5.10).  Let $x\in Drap(\cGx;\a_0)$
correspond to the filtration~(5.18). We argued that
$$\,\,\dim Drap(\cGx;\a_0)_x\le h^1(A,A)+h^1(B,A)+h^1(B,B)\,.\eqno(\star)$$
Here
$$\,\,h^1(B,B)=\dim Pic(S)+\dim Hilb^{\ell}(S)\,,$$
is the number of moduli of rank-one torsion-free sheaves with the
same Chern classes as $B$ ($\ell:=\ell(Z)=c_2(B)$).  For $m\in\RR$ let
$I_m(C)\ss Hilb^{\ell}(S)$ be the locus parametrizing subschemes $Z$ such
that $\sum_{P\in C}{\rm mult}_P(Z)\ge m$.  By~(5.21),
Inequality~($\star$) is replaced in the present case by
$$\dim Drap^C(\cGx;\a_0)_x\le
h^1(A,A)+h^1(B,A)+\dim Pic(S)+\dim I_{n\a\sqrt{H^2}}\,.$$
An easy application of a theorem of  Iarrobino~[I] shows
that
$$\,\,\cod\left(I_m(C),Hilb^{\ell}(S)\right)\ge m\,,$$
and hence
$$\dim Drap^C(\cGx;\a_0)_x\le
h^1(A,A)+h^1(B,A)+h^1(B,B)-n\a\sqrt{H^2}\,.$$
Since the right-hand side of the inequality in~(5.10) is an upper bound
for the above sum of $h^1$'s for $0\le\a\le\a_0$  we conclude
that~(5.11) holds.
\msk

\n
{\it Proof of Proposition~(5.15).}
\hskip 2mm
By additivity of the Euler characteristic, and by Serre duality, we have
$$\,\sum_{i\le j}h^1(Q_j,Q_i)=-\chi(F,F)+\sum_{i<j}\chi(Q_i,Q_j)+
\sum_{i\le j}h^0(Q_j,Q_i)+\sum_{i\le j}h^0(Q_i,Q_j\ot K_S)\,.$$
For $x\in {\rm Pic}(S)\ot\QQ$,  set
$\chi(x):=\chi(\cO_S)+(x^2-x\cdot K)/2$.
Applying the formula
$$\,\,\chi(F,G)=r_F\cdot r_G\left[\chi\left({c_1(G)\over r_G}-
{c_1(F)\over r_F}\right)-{\D_F\over r_F}-{\D_G\over r_G}\right]\,,$$
valid for any couple of sheaves $F$, $G$ (of positive rank) on $S$, we get:
$$\eqalign{\sum_{i\le j}h^1(Q_j,Q_i)=  & 2r_F\D_F-r_F^2\chi(\cO_S)+
\sum_{i<j}r_ir_j\left[\chi({c_1(Q_j)\over r_j}-
{c_1(Q_i)\over r_i})-{\D_i\over r_i}-{\D_j\over r_j}\right]\cr
&+\sum_{1<j}h^0(A,Q_j\ot K_S)+\sum_{2\le i\le j}h^0(Q_j,Q_i)\cr
&+h^0(A,A)+h^0(A,A\ot K_S)+\sum_{1<j}h^0(Q_j,A)
+\sum_{2\le i\le j}h^0(Q_i,Q_j\ot K_S)\,.}$$
Set
$$\eqalignno{\T:= &\sum_{i<j}r_ir_j\left[\chi\left({c_1(Q_j)\over r_j}-
{c_1(Q_i)\over r_i}\right)-{\D_i\over r_i}-{\D_j\over r_j}\right]\,, &
(5.22)\cr
\L:= & \sum_{1<j}h^0(A,Q_j\ot K_S)\,,\cr
\G:= & \sum_{2\le i\le j}h^0(Q_j,Q_i)\,, \cr
\O:= & h^0(A,A)+h^0(A,A\ot K_S)+\sum_{1<j}h^0(Q_j,A)
+\sum_{2\le i\le j}h^0(Q_i,Q_j\ot K_S)\,.}$$
Then the previous equation becomes
$$\,\,\sum_{i\le j}h^1(Q_j,Q_i)=2r_F\D_F-r_F^2\chi(\cO_S)+\T+\L+\G+\O\,.$$
Proposition~(5.15) follows from the bounds for $\T$, $\L$, $\G$ and $\O$
given below, together with a straightforward computation.

\proclaim (5.23).
Let notation be as above. If $K\cdot H \ge 0$, then
$$\eqalign{\T\le & -\D_F+{(r_A+1)r_B\over2r_A}\a^2
+\left[r_B(r_B-1)-r_F\right]{K\cdot H\over 2\sqrt{H^2}}\a \cr
& +\left[r_Fr_B(r_B-1)(r_F-2r_A)+r_Ar_B^2(r_A+1)\right]{(K\cdot
H)^2\over8H^2}\cr
& +M(r_F,S,H)+r_F^2\chi(\cO_S)\,,}$$

\proclaim (5.24).
With notation as above, assume $K\cdot H<0$. Then
$$\eqalign{\T\le & -\D_F+{(r_A+1)r_B\over2r_A}\a^2
-\left[(r_B-1)(2r_F-r_B)+r_F\right]{K\cdot H\over 2\sqrt{H^2}}\a\cr
 & +\left[r_Fr_B(r_B-1)(r_F+2r_A)+r_Ar_B^2(r_A+1)\right]
{(K\cdot H)^2\over8H^2}\cr
 & +{r_F^2\over 2}+{r_F^3\over 8}|K^2|\,.}$$

\n
In both cases (i.e.~regardless of the sign of $K\cdot H$), we have
$$\eqalignno{\L\le & {1\over 2r_A}(r_A^2+r_A+r_F)\a^2
+{r_F\over\sqrt{H^2}}(K\cdot H+s(r_F))\a+
{r_Ar_B\over2H^2}(K\cdot H+s(r_F))^2\,, & (5.25)\cr
\G\le & {r_B-1\over 2}\a^2+{r_B(r_B-1)s(r_F)\over \sqrt{H^2}}\a
+{r_B^3\over 2H^2}s(r_F)^2\,, &(5.26) \cr
\O\le & {r_Ar_F\over 2H^2}s(r_F)^2
+{r_A^2+r_B^2\over2H^2}(K\cdot H+s(r_F))^2\,. &(5.27)}$$
\bsk

\n
{\it Proof of~(5.23)-(5.24).}
\hskip 2mm
Let
$$\eqalignno{r_{\ell}:= &\max\{r_i\ |\ 1\le i\le n+1\}\,,\cr
\s_i:= & -r_1-\ldots-r_{i-1}+r_{i+1}+\ldots r_{n+1}\,. &(5.28)}$$
For future reference we notice that
$$\,\,\sum_{i=1}^{n+1}r_i\s_i=0\,.\eqno(5.29)$$

\proclaim (5.30) Lemma.
Keeping notation as above,
$$\,\,\T\le -(r_F-r_{\ell})\D_F
+{1\over 2}\left(r_F^2-\sum_{i=1}^{n+1}r_i^2\right)\chi(\cO_S)
+r_{\ell}\sum_{i=1}^{n+1}{c_1(Q_i)^2\over 2r_i}-r_\ell{c_1(F)^2\over 2r_F}
+{1\over 2}\sum_{i=1}^{n+1}\s_ic_1(Q_i)\cdot K_S\,.$$

\pf
We brake the sum defining $\T$ (Equation~(5.22)) into two pieces. First
consider
$$\,\,\T':=\sum_{i<j}r_ir_j\chi\left({c_1(Q_j)\over r_j}-{c_1(Q_i)\over
r_i}\right)\,.$$
A straightforward computation gives
$$\,\,\T'=r_F\sum_{i=1}^{n+1}{c_1(Q_i)^2\over 2r_i}-{1\over 2}c_1(F)^2+
{1\over 2}\sum_{i=1}^{n+1}\s_ic_1(Q_i)\cdot K_S+
{1\over 2}\left(r_F^2-\sum_{i=1}^{n+1}r_i^2\right)\chi(\cO_S)\,.$$
Now consider
$$\,\,\T'':=\sum_{i<j}(r_j\D_i+r_i\D_j)\,.$$
Then
$$\,\,\T''=\sum_{i=1}^{n+1}(r_F-r_i)\D_i\,.$$
Since $(r_F-r_i)\ge(r_F-r_{\ell})$, and since $\D_i\ge 0$ (Bogomolov's
theorem), we conclude that
$$\,\,\T''\ge (r_F-r_{\ell})\sum_{i=1}^{n+1}\D_i=
(r_F-r_{\ell})\sum_{i=1}^{n+1}[c_2(Q_i)-{1\over 2}c_1(Q_i)^2]+
(r_F-r_{\ell})\sum_{i=1}^{n+1}{1\over 2r_i}c_1(Q_i)^2\,.$$
Additivity of the Chern character  gives
$$\,\,\sum_{i=1}^{n+1}\left[c_2(Q_i)-{1\over 2}c_1(Q_i)^2\right]=
c_2(F)-{1\over 2}c_1(F)^2=\D_F-{1\over 2r_F}c_1(F)^2\,,$$
and thus we have
$$\,\,-\T''\le -(r_F-r_{\ell})\D_F+(r_F-r_{\ell}){1\over 2r_F}c_1(F)^2-
(r_F-r_{\ell})\sum_{i=1}^{n+1}{1\over 2r_i}c_1(Q_i)^2\,.$$
Since $\T=(\T'-\T'')$, the lemma follows.
\qed

\msk
\n
Set
$$\ \,\Xi:=r_{\ell}\sum_{i=1}^{n+1}{c_1(Q_i)^2\over 2r_i}-
r_\ell{c_1(F)^2\over 2r_F} +
{1\over 2}\sum_{i=1}^{n+1}\s_ic_1(Q_i)\cdot K_S\ .\eqno(5.31)$$

\proclaim (5.32) Lemma.
Let notation be as above. Then
$$\eqalign{\Xi\le &{r_{\ell}+r_Ar_{\ell}\over 2r_A}\a^2+
{1\over 2\sqrt{H^2}}\left\{\left[r_B(r_B-1)-r_F\right]
(K_S\cdot H)+r_F(r_B-1)\left[|K_S\cdot H|-(K_S\cdot H)\right]\right\}\a\cr
&+{r_B\over 8r_{\ell}H^2}\left\{\left[r_F^2(r_B-1)+r_Ar_B(r_A+1)\right]
(K_S\cdot H)^2-2r_Fr_A(r_B-1)|K_S\cdot H|(K_S\cdot H)\right\}\cr
&-{K^2\over 8r_{\ell}}\sum_{i=1}^{n+1}r_i\s_i^2\,.}$$

\pf
First of all notice that $\Xi$ is left invariant if $F$ and the $Q_i$ are
tensored
by a line bundle $\zeta$, or even if they are formally tensored by
$\zeta\in {\rm Pic}(S)\ot\QQ$. (Use Equation~(5.29).) Choosing $\zeta$
such that $c_1(F\ot\zeta)=0$, we can assume that $c_1(F)=0$. Of course by
doing this the classes $c_1(F)$ and $c_1(Q_i)$ become elements of
$NS(S)\ot\QQ$. Now rewrite the right-hand side of~(5.31) to get
$$\,\,\Xi={r_A\over 2r_{\ell}}\left[{r_{\ell}\over r_A}c_1(A)+{r_B\over
2}K\right]^2
+{1\over 2r_{\ell}}\sum_{i=2}^{n+1}
r_i\left[{r_{\ell}\over r_i}c_1(Q_i)+{\s_i\over 2}K\right]^2
-{K^2\over 8r_{\ell}}\sum_{i=1}^{n+1}r_i\s_i^2\,.\eqno(5.33)$$
In what follows we will use the inequality
$$\,\,L^2\le {(L\cdot H)^2\over H^2}\,,$$
which, by the Hodge index theorem, holds for all $L\in NS(S)\ot\QQ$. It   gives
$$\,\,{r_A\over 2r_{\ell}}\left[{r_{\ell}\over r_A}c_1(A)+{r_B\over
2}K\right]^2
\le{r_{\ell}\over 2r_A}\a^2-{r_B(K\cdot H)\over2\sqrt{H^2}}\a+
{r_Ar_B^2(K\cdot H)^2\over 8r_{\ell}H^2}\,,\eqno(5.34)$$
and
$$\,\,{1\over 2r_{\ell}}\sum_{i=2}^{n+1}
r_i\left[{r_{\ell}\over r_i}c_1(Q_i)+{\s_i\over 2}K\right]^2
\le {1\over 2r_{\ell}H^2}\sum_{i=2}^{n+1}
r_i\left[r_{\ell}\mu_i+{\s_i\over 2}\mu_K\right]^2\,.\eqno(5.35)$$
We will bound the right-hand side of the above inequality by applying the
following lemma. Its (easy) proof is left to the reader.

\proclaim (5.36) Lemma.
Let $x_1,\ldots,x_n$ be real numbers and $r_1,\ldots,r_n$ be
positive integers. Let $N:=\sum_{i=1}^nr_i$. If $x_i\ge a$ for
$i=1,\ldots,n$, then
$$\,\,\sum_{i=1}^nr_ix_i^2\le[\sum_{i=1}^nr_ix_i-(N-1)a]^2+(N-1)a^2\,.$$

\bsk
\n
We apply the lemma to the sum on the right-hand side of~(5.35).
Set   $x_i:=r_{\ell}\mu_{i+1}+\s_{i+1}\mu_K/2$. Then
$N=r_B$,  and
$$\,\,\sum_{i=2}^{n+1}r_ix_i=
r_{\ell}\sqrt{H^2}\a-{1\over 2}r_Ar_B\mu_K\,.$$
Since $\mu_2\ge\cdots\ge\mu_{n+1}\ge 0$,  we can  set
$a:=-{1\over 2}r_F|\mu_K|$.  Lemma~(5.36) gives
$$\eqalign{{1\over 2r_{\ell}H^2}\sum_{i=2}^{n+1}
r_i\left[r_{\ell}\mu_i+{\s_i\over 2}\mu_K\right]^2\le & {1\over 2}r_{\ell}\a^2+
{1\over 2\sqrt{H^2}}\left[r_F(r_B-1)|\mu_K|-r_Ar_B\mu_K\right]\a\cr
& +{1\over 8r_{\ell}H^2}\left\{\left[r_F(r_B-1)|\mu_K|-r_Ar_B\mu_K\right]^2
+r_F^2(r_B-1)\mu_K^2\right\}\,.}$$
Lemma~(5.32) now follows from~(5.33), (5.34), (5.35), and the above inequality,
together with   a straightforward computation.
\qed

\msk
Inequalities~(5.23)-(5.24) follow from Lemmas~(5.30)-(5.32) and
some easy estimates. The only estimate which is not completely trivial is
provided by the following
\bsk

\proclaim Lemma.
Keeping notation as above, we have
$$\,\,r_F\le {1\over r_{\ell}}\sum_{i=1}^{n+1}r_i\s_i^2\,.$$

\pf
One checks easily that, if $n$ is replaced by $(n+1)$, and $r_{n+1}$  by
$t_{n+1},t_{n+2}$ (with $t_{n+1}+t_{n+2}=r_{n+1}$), then the right-hand side of
the
above inequality increases. Thus the minimum, with a fixed $r_1=r_A$, is given
by
$n=1$. The lemma then follows by direct computation.
\qed
\msk

\n
{\it Proof of Inequality~(5.25).}
\hskip 2mm
By Corollary~(5.3), we have
$$\,\,\L\le{r_A\over 2H^2}
\sum_{1<j}r_j(\mu_j+\mu_K-\mu_A+s(r_F))^2\,.$$
Lemma~(5.36) together with an easy computation gives~(5.25).
\msk

\n
{\it Proof of Inequality~(5.26).}
\hskip 2mm
To simplify notation   let $s:=s(r_F)$.   By Corollary~(5.3) we have
$$\,\,\G\le {1\over 2H^2}\sum_{1<i<j}r_ir_j\left(\mu_i-\mu_j+s\right)^2
+\sum_{1<i}r_i^2s^2\,.$$
Expanding the squares in the first sum on the right-hand side, we write
$$\,\,\G={1\over 2H^2}\left(\G_1+\G_2+\G_3\right)\,,\eqno(5.37)$$
where
$$\eqalign{\G_1:= & s^2\sum_{1<i}r_i^2+s^2\sum_{1<i<j}r_ir_j\,,\cr
 \G_2:= & \sum_{1<i<j}\left(r_ir_j\mu_i^2+r_ir_j\mu_j^2\right)
-2\sum_{1<i<j}r_ir_j\mu_i\mu_j\,,\cr
\G_3:= & 2s\sum_{1<i<j}\left(r_ir_j\mu_i-r_ir_j\mu_j\right)\,.}$$
We rewrite $\G_2$ as
$$\,\,\G_2=\sum_{1<i}r_i(r_B-r_i)\mu_i^2-2\sum_{1<i<j}r_ir_j\mu_i\mu_j
=r_B\sum_{1<i}r_i\mu_i^2-r_B^2\mu_B^2\,.$$
For the second equality we have used the relation
$r_B\mu_B=\sum_{1<i}r_i\mu_i$.
We also write
$$\,\,\G_3=2s\sum_{1<i}r_i\tau_i\mu_i\,,$$
where $\tau_i:=(-r_2-\cdots-r_{i-1}+r_{i+1}+\cdots+r_{n+1})$. (See~(5.28).) Now
write
$$\eqalign{\G_2+\G_3 = &
r_B\sum_{1<i}r_i\left[\left(\mu_i+{\tau_is\over r_B}\right)^2
-{\tau_i^2s^2\over r_B^2}\right] -r_B^2\mu_B^2 \cr
\le & r_B\sum_{1<i}r_i\left(\mu_i+{\tau_is\over r_B}\right)^2
-r_B^2\mu_B^2\,.}\eqno(5.38)$$
We bound the sum of squares on the last line by applying Lemma~(5.36).
We let $a:=(\mu_F-s)$. Clearly $N=r_B$. Furthermore
$$\,\,\sum_{1<i}r_i\left(\mu_i+{\tau_is\over r_B}\right)
=r_B\mu_B+{s\over r_B}\sum_{1<i}r_i\tau_i=r_B\mu_B\,.$$
Lemma~(5.36)  gives
$$\,\,\sum_{1<i}r_i\left(\mu_i+{\tau_is\over r_B}\right)^2
\le r_B^2\mu_B^2-2r_B(r_B-1)\mu_B(\mu_F-s)+r_B(r_B-1)(\mu_F-s)^2\,.$$
Putting this together with Inequality~(5.38) we conclude that
$$\,\,\G_2+\G_3\le r_B^2(r_B-1)(\mu_B-\mu_F+s)^2
=(r_B-1)H^2\a^2+2r_B(r_B-1)s\sqrt{H^2}\a+r_B^2(r_B-1)s^2\,.$$
Adding $\G_1\le r_B^2s^2$, and using Equation~(5.37), one gets~(5.26).
\msk

\n
{\it Proof of Inequality~(5.27).}
\hskip 2mm
It follows from Corollary~(5.3) and some simple estimates.
\msk

\n
{\it Proof of Proposition~(5.16).}
\hskip 2mm
The proof follows that of Proposition~(5.15), with the difference that we
replace Corollary~(5.3) by Proposition~(5.5).
\bsk

\n
{\bf Properly semistable sheaves.}
\msk
\n
Let $F$ be  a torsion-free sheaf on $S$ which is properly $\mu$-semistable,
i.e. $\mu$-semistable but not $\mu$-stable.    We
will bound the dimension of the locus $V^0(F)\ss Def^0(F)$ parametrizing
properly $\mu$-semistable sheaves.

\proclaim (5.39).
For $r\ge 2$  an integer, set
$$\e(r,S,H):=\cases{{r^2\over 2H^2}\left[K\cdot H+s(r)\right]^2
+{r^2\over 16}\left[{(K\cdot H)^2\over H^2}-K^2\right]
+r^2+pr^2|\chi(\cO_S)|-q_S
& if $r>2$, \cr
{3\over 2H^2}\left(K\cdot H+H^2+1\right)^2
+{1\over 8}\left[{(K\cdot H)^2\over H^2}-K^2\right]
+4-3\chi(\cO_S)-q_S & if $r=2$,}$$
where $p:=-3/4$ if $\chi(\cO_S)\ge 0$, and $p:=1$ if $\chi(\cO_S)<0$.

\proclaim (5.40) Proposition.
Let $(S,H)$ be a polarized surface. Let $F$ be a torsion-free sheaf on $S$ of
rank $r_F\ge 2$, which is $\mu$-semistable but not $\mu$-stable
(i.e.~properly $\mu$-semistable). If $r_F>2$ assume that $H$
satisfies~(0.4), if $r_F=2$ assume only that $H$ is effective. Then,
letting $V^0(F)\ss Def^0(F)$ be as above, we have
$$\,\,\dim V^0(F)\le \left(2r_F-1\right)\D_F+\e(r_F,S,H)\,.$$

\n
Notice that if $F$ is simple, then the proof that Proposition~(5.15)
implies Proposition~(5.9) gives also a bound for $\dim V^0(F)$.  If however
$F$ is not simple, then we need to argue  differently. Before proving the
above bound, we state a lemma.

\proclaim (5.41) Lemma.
Let $A$, $B$ be torsion-free sheaves on a projective irreducible variety.
Suppose that
$A$, $B$ are both properly $\mu$-semistable, and that $\mu_A=\mu_B$. Then
$$\,\,h^0(A,B)\le r_Ar_B\,.$$

\pf
The proof is by double induction on the lengths of  $\mu$-Jordan-H\"older
filtrations
for $A$ and $B$. If the lengths are both one, then $A$ and $B$ are stable,
hence
$h^0(A,B)\le 1$, and the result is true. To prove the inductive step, let
$A_1\ss A$
and $B_1\ss B$ be the first terms of Jordan-H\"older filtrations. Then
$$\,\,h^0(A,B)\le 1+r_{A_1}(r_B-r_{B_1})+r_{B_1}(r_A-r_{A_1})
+(r_A-r_{A_1})(r_B-r_{B_1})\,.$$
The lemma follows by simplifying the right-hand side.
\qed
\msk

\n
{\it Proof of Proposition~(5.40).}
\hskip 2mm
Let  $V(F)\ss Def(F)$ be the locus parametrizing properly
$\mu$-semistable sheaves. We will prove that
$$\,\,\dim V(F)\le (2r_F-1)+\e(r_F,S,H)+q_S\,.\eqno(5.42)$$
Clearly this is equivalent to Proposition~(5.40). Let $\cF$ be the family of
sheaves
parametrized by $Def(F)$. Let $Quot_0(\cF)$ be the quot-scheme parametrizing
torsion-free quotients of $\cF_x$, for $x$ varying in $Def(F)$, with slope
equal to
that of $F$ (i.e.~destabilizing quotients).  By a theorem of Grothendieck~[G],
$Quot_0(\cF)$ is of finite type.  Let
$$\,\,\pi\cl Quot_0(\cF)\to Def(F)\,,$$
be the projection. Clearly
$$\,\,V(F)=\pi\left(Quot_0(\cF)\right)\,,\eqno(5.43)$$
and hence $\dim V(F)\le \dim Quot_0(\cF)$. Now let $y\in Quot_0(\cF;P)$
correspond to the $\mu$-detabilizing sequence
$$\,\,0\to A\to \cF_x\to B\to 0\,.\eqno(*)$$
Then there is an exact sequence~[DL]
$$\,\,0\to Hom(A,B)\to T_yQuot_0(\cF_x;P)\brel \pi_*\over\lra
Ext^1(\cF_x,\cF_x)
\brel\o_{+}\over\lra Ext^1(A,B)\,,\eqno(5.44)$$
where $\o_{+}$ is induced by the inclusion $A\ss \cF_x$ and the quotient
$\cF_x\to B$.
By applying the functors $Hom(\cdot,\cF_x)$ and $Hom(A,\cdot)$ to~($*$) one
concludes
that  %
$$\,\,\dim Quot_0(\cF_x;P)_y\le h^0(A,B)+h^1(A,A)+h^1(B,B)+h^1(B,A)\,.$$
Now we estimate the sum of the $h^1$'s by repeating the proof of
Proposition~(5.15).
(We will get sharper results because the filtration consists of only two
terms.) Adopting the notation used in that proof, we have
$$\,\,\T\le  -r_A\D_F+{r_Ar_F\over 8}\left[{(K\cdot H)^2\over
H^2}-K^2\right]\,.$$
Then, by applying~(5.3), (5.5), or~(5.41) to bound $\L$, $\G$ and
$\O$, we get
$$\,\,\dim Quot_0(\cF;P)_y\le (2r_F-1)+\e(r_F,S,H)+q_S\,.$$
By~(5.43) this proves Inequality~(5.42), and hence the proposition.

\proclaim (5.45) Corollary.
Let $(S,H)$ be a polarized surface, and  $\xi$ be  set of sheaf data for $S$.
 If $\rx>2$ we assume that $H$ satisfes~(0.4), if $\rx=2$ we only
assume that $H$ is effective.  Let
$X\ss\cMx$ be a subset such that
$$\,\,\dim X>(2\rx-1)\Dx+\e(\rx,S,H)\,.$$
Then there exists a point of $X$ parametrizing a $\mu$-stable sheaf.

\pf
Let $[F]\in X$; we can assume $F$ is properly $\mu$-semistable. Let  $\cF$ be
the family of sheaves on $S$ parametrized by $Def^0(Gr F)$. By Luna's \`etale
slice Theorem the map
$$\l\cl Def^0(Gr F)\to \cMx\,,$$
induced by $\cF$ is surjective onto a neighborhood of $[F]$. Thus
$\dim\l^{-1}X$ satisfies
the same inequality as $\dim X$. By~(5.40) there exists $x\in\l^{-1}X$
parametrizing a
$\mu$-stable sheaf. Then $\l(x)\in X$ is a point parametrizing a $\mu$-stable
sheaf.
\qed

\proclaim (5.46) Corollary.
Let $(S,H)$ and $\xi$ be as in the previous corollary. If
$$\,\,2\rx\Dx-(\rx^2-1)\chi(\cO_S)>(2\rx-1)\Dx+\e(\rx,S,H)\,,$$
then the generic point of any irreducible component of $\cMx$ parametrizes
a $\mu$-stable sheaf.

\pf
Let $X\ss\cMx$ be an irreducible component. Let $[F]\in X$ be a point not
belonging to any other component of $\cMx$.  If $F$ is $\mu$-stable  there
is nothing to prove, so assume $F$ is not $\mu$-stable.   By deformation
theory~[F] we have
$$\,\,\dim Def^0(F)\ge 4\Dx-3\chi(\cO_S)\,.$$
Hence by Proposition~(5.40) we conclude that the generic
point  $x\in Def^0(F)$ parametrizes a $\mu$-stable sheaf. This implies
that the generic point of $X$ parametrizes a $\mu$-stable sheaf.
\qed
\bsk

\n
{\bf Non-stable vector bundles on curves.}
\msk
\n
Let $C$ be a smooth irreducible curve of genus $g$. Let $\cF$ be a family of
rank-$r$
vector bundles on $C$, parametrized by an equidimensional variety $B$. Let
$B^{ns}\ss B$ be the subset parametizing bundles which are not stable
(i.e.~either unstable or properly semistable). The goal of this subsection is
to prove the following

\proclaim (5.47) Proposition.
Keeping notation as above, assume that $B^{ns}$ is not empty.
Then
$$\,\,\cod (B^{ns},B)\le {r^2\over 4}g\,.$$

\n
If $g=0$, there are no stable vector bundles, and hence the proposition is
trivially verified. Thus we can assume that $g>0$.
Let $F$ be a non-stable vector bundle on $C$. Let $V(F)\ss Def(F)$ be the
subset parametrizing non-stable bundles. Since $Def(F)$ is smooth, it
suffices  to show that, if $V_i$ is an  irreducible component of $V(F)$, then
$$\,\,\cod \left(V_i,Def(F)\right)\le {r^2\over 4}g\,.\eqno(5.48)$$
This is what we will  prove. One can stratify $V(F)$
according to the ranks and slopes of the successive quotients of the
Harder-Narasimhan
filtration. Thus the stratum corresponding to the {\it type}
$$\,\,{\bf t}:=\left((r_1,\mu_1),\ldots,(r_n,\mu_n)\right)\,,$$
where $\mu_1>\cdots >\mu_n$, consists of the points $x\in Def(F)$ such
that there is a filtration
$$\,\,F_1\ss F_2\ss\cdots\ss F_n=\cF_x\,,$$
with $F_i/F_{i-1}$ a rank-$r_i$ semistable bundle with slope $\mu_i$. As is
well-known~[AB,Le], each $V_{\bf t}$ is smooth, and
$$\,\,\cod\left(V_{\bf t},Def(F)\right)
=\sum_{i<j}r_ir_j(\mu_i-\mu_j+g-1)\,.\eqno(5.49)$$
In order to prove Inequality~(5.48) we need the following

\proclaim  Lemma.
Keep notation as above. Assume that the genus $g$ of $C$ is positive.
Let $F$ be a non-stable vector bundle on $C$. Then every irreducible
component of $V(F)$ contains an open dense subset parametrizing minimally
non-stable bundles, i.e.~bundles $F$ fitting into an exact sequence
$$\,\,0\to A\to F\to B\to 0\,,\eqno(5.50)$$
where $A$, $B$ are semistable vector bundles such that
$$\eqalignno{\mu_A & \ge \mu_B & (5.51)\,,\cr
\mu_A-{1\over r_A} & <  \mu_B+{1\over r_B} & (5.52)\,.}$$

\pf
Fix a component $V_i$ of $V(F)$. If the  generic point of
$V_i$ parametrizes a  semistable bundle, then there is nothing to prove. Thus
we can
assume that the bundles parametrized by $V_i$ are unstable. Let ${\bf t}$ be
the type of
the Harder-Narasimhan filtration for the generic bundle parametrized by $V_i$.
Since versality is an open condition, we can  replace $F$ by $\cF_x$ for  a
generic point $x\in V_i$,  and thus we can assume that $V(F)$ is irreducible
and that $V(F)=V_{\bf t}$. First let's show that the Harder-Narasimhan
filtration corresponding to ${\bf t}$ consists of only two terms. Let $A$ be
the first term of the H.-N. filtration for $F$, and let
$$\,\,0\to A\to F\to B\to 0\,,\eqno(*)$$
be the corresponding exact sequence. Let $P$ be the Hilbert polynomial of
$B$, and let $Quot(\cF;P)$ be the $Quot$-scheme parametrizing quotients of
$\cF_x$, for
$x\in Def(F)$, with Hilbert polynomial equal to $P$. If $y\in Quot(\cF;P)$,
then there is
an exact sequence~(5.44). The last map is surjective, because $C$ is smooth of
dimension one. By a criterion of Drezet-LePotier~[DL], $Quot(\cF;P)$ is
smooth. Hence its dimension can be computed from~(5.44). Furthermore,
since $h^0(A,B)=0$, the projection $\pi\cl Quot(\cF;P)\to Def(F)$ is an
embedding; let $W$ be its image.  An easy computation gives
$$\,\,\cod\left(W,Def(F)\right)=r_Ar_B(\mu_A-\mu_B+g-1)\,.$$
Clearly $W\ss V(F)$, and hence
$\cod\left(W,Def(F)\right)\ge\left(V(F),Def(F)\right)$. Since $V(F)=V_{\bf
t}$,  one easily concludes from Equation~(5.49) that ${\bf
t}=\left((r_A,\mu_A),(r_B,\mu_B)\right)$. Now  let's show that~($*$) is
minimally destabilizing, i.e.~that~(5.52) is satisfied. We argue by
contradiction. Assume that~(5.52) is violated. Let $A'$ be a sheaf fitting
into an exact sequence
$$\,\,0\to A'\to A\brel\phi\over\lra k_P\to 0\,,\eqno(5.53)$$
where $k_P$ is the skyscraper sheaf at  some point $P\in C$. Let $B':=B\op
k_P$.
Then there is an exact sequence
$$\,\,0\to A'\to F\to B'\to 0\,.\eqno(\dag)$$
Let $P'$ be the Hilbert polynomial of $B'$, and let $Quot(\cF;P')$ be the
$Quot$-scheme
of quotients of $\cF_x$ with Hilbert polynomial equal to $P'$, for $x\in
Def(F)$. As in
the previous case, $Quot(\cF;P')$ is smooth, and~(5.44) gives
$$\,\,\dim
Quot(\cF;P')=h^1(F,F)-r_Ar_B(\mu_A-\mu_B+g-1)+r_F\,.\eqno(\sharp)$$
Let $\pi\cl Quot(\cF;P')\to Def(F)$ be the projection.  Since~(5.52)  is
violated,
$$\,\,\pi\left(Quot(\cF;P')\right)\ss V(F)\,.\eqno(\star)$$
 As is easily checked,
$\pi^{-1}({\bf o})$ consists of all sequences~($\dag$), where $A'$ fits into an
exact
sequence~(5.53). (Here ${\bf o}$ is the point parametrizing $F$.) Letting  $P$
and $\phi$ vary, we get
$$\,\,\dim\pi^{-1}({\bf o})=r_A\,.$$
Since $r_A<r_F$,  we conclude by~($\sharp$) that
$$\,\,\dim \pi\left(Quot(\cF;P')\right)>
h^1(F,F)-r_Ar_B(\mu_A-\mu_B+g-1)\,.$$
This inequality, together with~($\star$), contradicts Formula~(5.49). This
proves that~(5.52) is satisfied, and hence $F$ is minimally non-stable.
\qed
\msk

\n
{\it Proof of~(5.48).}
\hskip 2mm
By the previous lemma we can assume that every sheaf parametrized by
$V_i$ fits into  Exact sequence~(5.50), with~(5.51)-(5.52)  satisfied. We
distinguish two cases, according to whether $V_i$  parametrizes unstable or
semistable sheaves. In the first case Formula~(5.49) gives
$$\cod\left(V_i,Def(F)\right)=r_Ar_B(\mu_A-\mu_B+g-1)\,.$$
By~(5.51)-(5.52) we conclude that~(5.48) is satisfied. Now assume that $F$ is
semistable. Let $Quot_0(\cF)$ be the $Quot$-scheme parametrizing quotients of
$\cF_x$
whose slope is equal to that of $\cF_x$, for $x\in Def(F)$. Let $\pi\cl
Quot_0(\cF)\to
Def(F)$ be the projection.  Exact sequence~(5.44) and the smoothness of
$Quot_0(\cF)$
give
$$\,\,\dim \pi\left(Quot_0(\cF)\right)\ge h^1(F,F)+\chi(A,B)-h^0(A,B)\,.$$
Since $\pi\left(Quot_0(\cF)\right)=V_i$ one concludes, by Lemma~(5.41),
that~(5.48) is satisfied.
\bsk
\bsk

\n
{\lbf References.}
\msk

\frenchspacing
\item{[AB]} M. F. Atiyah and R. Bott. The Yang-Mills equations over Riemann
surfaces. {\it Philosophical Transactions of the Royal Society of London,
Series A}, 308, 523-615.
\item{[BPV]} W. Barth, C. Peters and A. Van de Ven. {\it Compact complex
surfaces}, Ergebnisse der Mathematik und ihrer Grenzgebiete 3. Folge-Band
4 (1984).
\item{[Bo]} E. Bombieri. Canonical models of surfaces of general type, {\it
Publ. Math. Inst. Hautes Etud. Sci.} 42 (1973), 171-219.
\item{[D]} S. K. Donaldson. Polynomial invariants for smooth
four-manifolds,  {\it Topology} 29 (1990), 257-315.
\item{[DL]} J. M. Drezet and J. Le Potier. Fibr\'es stables et fibr\'es
exceptionnels sur le plan projectif,{\it Ann. scient. Ec. Norm. Sup.} ${\rm
4}^e$ s\'erie t. 18 (1985), 193-244.
\item{[DN]} J. M. Drezet and M. S. Narasimhan. Groupe de Picard des
vari\'et\'es de modules de fibr\'es semi-stables sur les courbes
alg\'ebriques, {\it Invent. math.} 97 (1989), 53-94.
\item{[G1]} D. Gieseker. On the moduli of vector bundles on an algebraic
surface, {\it Ann. of Math.} 106 (1977), 45-60.
\item{[G2]} D. Gieseker. A degeneration of the moduli space of stable
bundles, {\it J. Differ. Geom.} 19 (1984), 173-206.
\item{[GL1]} D. Gieseker and J. Li. Irreducibility of moduli of rank two
vector bundles, to appear in {\it J. Differ. Geom.}
\item{[GL2]} D. Gieseker and J. Li.  Moduli of vector bundles over surfaces I,
preprint.
\item{[G]} A. Grothendieck, Techniques de construction et th\'eore\`emes
d'existence en g\'eometrie algebrique, IV: les sch\'emas de Hilbert, {\it
S\'em. Bourbaki} 221 (1960).
\item{[HL]} A. Hirschowitz and Y. Laszlo. A propos de l'existence de fibr\'es
stables sur les surfaces, {\it preprint}.
\item{[I]} A. Iarrobino.  Punctual Hilbert schemes, {\it  Bulletin of the
American Mathematical Society} 78 (1972), 819-823.
\item{[Le]} J. Le Potier. Espaces de modules de faisceaux semi-stables sur
le plan projectif, {\it preprint} School "Vector bundles on surfaces"-CIMI
and Europroj, Nice-Sophia-Antipolis June 1993.
\item{[Li]} J. Li. Algebraic geometric interpretation of Donaldson's
polynomial invariants, {\it J. Differ. Geom.} 37 (1993) 417-466.
\item{[LQ]} W. P. Li and Z. Qin. Stable vector bundles on algebraic surfaces,
{\it preprint}.
\item{[Ma]} M. Maruyama. Moduli of stable sheaves II, {\it J. Math. Kyoto
Univ.} 18-3 (1978), 557-614.
\item{[MO]} J. Morgan and K. G. O'Grady. Differential topology of complex
surfaces Elliptic surfaces with $p_g=1$: smooth classification,  {\it
Lecture Notes in Mathematics} 1545, Springer-Verlag.
\item{[Mu]} S. Mukai. On the moduli space of bundles on $K3$ surfaces I,
Vector bundles on algebraic varieties, {\it Tata Institute of fundamental
research studies in mathematics}, Oxford University Press (1987).
\item{[O]} K. G. O'Grady. The irreducible components of moduli spaces of
vector bundles on surfaces, {\it Invent. Math.} 112 (1993), 585-613.
\item{[S]} C. Simpson. Moduli of representations of the fundamental group
of a smooth projective variety I, {\it preprint}.
\item{[Zuo]} K. Zuo. Generic smoothness of the moduli of rank two stable
bundles over an algebraic surface, {\it preprint MPI/90-7}.

\end